\documentclass[11
pt]{article}

\newcounter{alfej}[section]

\newcommand{\bekezdes}[2]{\vspace*{-2mm}\paragraph{#1}\label{#2}}

\usepackage{amssymb}

\begin{document}
\let\n\noindent
\let\d\partial
\def\EE{\mathcal E}
\def\C{\mathbb C}
\def\Q{\mathbb Q}
\def\R{\mathbb R}
\def\Z{\mathbb Z}
\def\N{\mathbb N}
\def\Res{\operatorname{Res}}
\def\im{\operatorname{im}}
\def\Ker{\operatorname{Ker}}
\let\ker\Ker
\def\Tot{\operatorname{Tot}}
\def\Gr{\operatorname{Gr}}
\def\Hom{\operatorname{Hom}}
\def\Log{\operatorname{Hom}}
\def\codim{\operatorname{codim}}
\let\emptyset\varnothing
\def\resp{\operatorname{respectively}}
\def\resp{\operatorname{respectively}}
\def\mod{\operatorname{mod}}
\def\modulo{\operatorname{modulo}}
\def\b{\bullet}
\def\c{\mathcal}
\def\p{\pitchfork}
\def\Y{Y_{\alpha_1,\ldots,\alpha_p}}
\def\u{\underline}
\def\y{\tilde{Y}^p}
\def\s{\star}
\def\noi{\noindent}
\newcommand{\be}{\bekezdes}
\def\esq{$\hfill\square$}
%line 53

\begin{center}
{\Large \bf Deligne-Hodge-DeRham theory with coefficients}\\
 \ \\ {\Large Fouad Elzein}
\end{center}

\begin{abstract}
 Let ${\cal L}$ be a variation of Hodge structures  on the complement $X^{*}$ of
 a normal crossing divisor (NCD) $ Y$ in a smooth analytic variety $X$
and let $ j: X^{*} = X - Y \to X $ denotes the open embedding. The
purpose of this paper is to describe  the weight filtration $W$ on
a combinatorial logarithmic
 complex computing the (higher)
 direct image ${\bf j}_{*}{\cal L} $, underlying a mixed Hodge complex when $X$ is
proper, proving in this way the results in the note [14]
generalizing the constant coefficients case. When a morphism $f: X
\to D$ to a complex disc is given  with $Y = f^{-1}(0)$,
 the weight filtration on
the complex of nearby cocycles $\Psi_f ({\cal L})$ on $Y$ can  be
described by these logarithmic techniques and a comparison theorem
shows that the filtration coincides with the weight defined by the
logarithm of the monodromy which provides the link with various
results on the subject.

\vspace{2mm}
\n { \bf Mathematics Subject Classification(2000)}:
32S35 (primary), 32G20, 14D07 (secondary).

\end{abstract}
% paragraphe 2 : \S 2
\section{Introduction}
The subject of this article is to construct a mixed Hodge
structure $MHS$ on the cohomology of a  local system ${\cal L}$
underlying a polarised variation of Hodge structures ($VHS$) on
the complement $X^{*}= X - Y$ of a normal crossing divisor (NCD) $
Y$ in a smooth proper complex algbraic or analytic variety $X$.
 Let $ j: X^{*}  \to X $ denotes
the open embedding.
 Technically, we need to define a structure of mixed Hodge complex
 $MHC$ on  the higher direct image ${\bf j}_{*}{\cal L}$.\\
 The work consists then in two parts,
  first to define the
rational weight filtration
 $ W $
 and second to construct the complex weight $ W $ and Hodge filtration $ F $.
 Although we will use the same letter ${\cal L}$ for the rational
 as well the complex local system, sometimes when we need to stress the difference
 we denote by ${\cal L}^r$ the rational and by ${\cal L}^c$ the
 complex local system. \\
 In the rational case, we don't have a particular representative
 of ${\bf j}_{*}{\cal L}$ by a distinguished complex, so the method
 is to use the theory of perverse sheaves to describe $W$.\\
 While in the complex case we need to construct a bi-filtered
 complex and we use the logarithmic complex with coefficients in
 {\it Deligne's analytic extension with regular singular
connection ${\cal L}_X $}
the bundle extension of  $ {\cal L} \otimes {\cal O}_{X^*} $ since
by Deligne's theorem
$${\bf j}_{*}{\cal L} \cong
  {\Omega}^*_X (Log Y) \otimes {\cal L}_X \leqno{(1)}$$
\n By the subsequent work of Schmid, Cattani and Kaplan, Kashiwara
and Kawai, the   Hodge filtration $F$ extends by sub-bundles.
 \n In this article we describe a bi-filtered complex
$$ ({\Omega}^* {\cal L},{\cal W} ,F )\leqno{(2)}$$
\n constructed as a sum of a combinatorial complex constantly
equal to
 the logarithmic complex with coefficients and
  which underly the structure of $MHC$ we are looking for.
  Although the existence of such bi-filtered complex is important
  in the general theory, the basic results on $ Gr^{\cal W} {\cal L}$ can be stated
 more easily so to reflect the topological and geometrical
 properties of the variety and the local system.
  Let us fix the hypothesis and the notations for the rest of
 the article.

\smallskip
\n {\bf  Hypothesis }. {\it  Let  ${\cal L}$ be a polarised  local
system defined over $\Q$, on the complement of the normal crossing
divisor (NCD) $ Y$ in a smooth analytic variety $X$, $({\cal L}_X,
\nabla) $ the canonical extension  of ${\cal L}_{X^*} = {\cal L}
 \otimes {\cal O}_{X^*}  $ [6] with a meromorphic connection $\nabla$
 on $X$ having a regular singularity along $ Y$ in $X$ and
 $ {\Omega}^*_X (Log Y) \otimes {\cal L}_X $ the associated DeRham
 logarithmic complex  defined by
$\nabla$, moreover  ${\cal L}^{\C} = {\cal L} \otimes \C $
underlies a polarised variation of Hodge structures
($VHS$).}\\
In the text  we write ${\cal L}$ for the rational ${\cal L}^r$ as
well complex local system ${\cal L}^c$, ${\cal L}_X$ for the
analytic extension  and to simplify the exposition,
 we suppose ${\cal L}$ locally unipotent along $Y$.

\smallskip
\n {\bf Notation}. {\it We suppose the NCD, $Y = \cup_{i \in I}
Y_i $ equal to the union of irreducible and smooth components
$Y_i$ for $i$ in $I$. For all subset $K$ of $I$, let $Y_K =
{\cap_{i \in K}} Y_i$, $Y_K^* = Y_K - {\cup_{i \not\in K}} Y_i$,
and $j^K \colon Y^*_K \rightarrow Y_K $  the locally closed
embedding, then $Y_K - Y_K^*$ is a NCD in $Y_K$ and the open
subsets $Y_K^*$ of $Y_K$ form with $X^*$ a natural stratification
of $X$. All extensions of $\cal L$ considered are constructible
with respect to this
stratification and even perverse.\\
We write ${\cal L}_{Y_K}$ for the restriction of $ {\cal L}_X$ to
$Y_K$, ${\cal N}_i$ for the nilpotent endomorphisms of the
restriction ${\cal L}_{Y_K}$ logarithm of the unipotent part of
the monodromy, and
$$ {\cal W}_{Y_K}^K = {\cal W}(\Sigma_{i\in K} {\cal
N}_i)\leqno{(3)}$$
  \n for the filtration by sub-bundles defined on
${\cal L}_{Y_K}$ by $\Sigma_{i\in K} {\cal N}_i$.}

\smallskip
\n   Let $i_K:Y^{*}_K \rightarrow X $, we introduce the local
systems

\smallskip
\centerline { ${\cal L}^K  = i_K^* R^{\vert K \vert}{\bf j}_*
 {\cal L}$, \quad
  ${\cal L}'^K = i_K^* R^0 {\bf j}_* {\cal L}$ }

\n and the bundles

\smallskip
\centerline {  ${\cal L}^K_{X} =
 {\cal L}_{Y_K}/(\Sigma_{i \in K} {\cal N}_i {\cal L}_{Y_K}) $,
 \quad
${\cal L}'^K_{X} = \cap_{i \in K} (ker {\cal N}_i: {\cal L}_{Y_K}
\to {\cal L}_{Y_K}) $ }

\smallskip
\n Throughout this work we prove the following results

\medskip
 \n {\bf Theorem 1}. {\em  i) ${\cal L}^K_{X} $  (resp. ${\cal L}'^K_{X}$)
 is a flat bundle with flat sections isomorphic to  the  local system
 ${\cal L}^K$ (resp. ${\cal L}'^K$); precisely  they are resp. Deligne's
extension of the complex local system.

\smallskip
\n ii) The filtration $ {\cal W}_{Y_K}^K $ (3) induces a filtration by
flat sub-bundles of ${\cal L}^K_{X} $ (resp ${\cal L}'^K_{X} $),
hence induces a filtration by complex sub-local systems $ {\cal
W}^K $ of ${\cal L}^K$ ( resp. ${\cal L}'^K $).

\smallskip
\n iii) The filtration $ {\cal W}^K $ is defined on the rational
local system $ {\cal L}^K $ ( resp. ${\cal L}'^K $).

\smallskip
\n iv)  Let ${\cal L}_r^K \colon= Gr^{{\cal W}^K}_{r - \vert K
\vert }{\cal L}^K$ and ${\cal L}_{X,r}^K \colon= Gr^{{\cal
W}_{Y_K}^K}_{ r -
\vert K \vert }{\cal L}^K_{X}$ for $r>0$, resp.\\
 ${\cal L}_r^K \colon= Gr^{{\cal W}^K}_{r + \vert K \vert }{\cal L}'^K$
 and ${\cal L}_{X,r}^K \colon=
Gr^{{\cal W}_{Y_K}^K}_{ r + \vert K \vert }{\cal L}'^K_{X}$ for
$r<
0$,\\
 then the  system defined by  $({\cal L}_r^K, {\cal L}_{X,r}^K , F)$
 where $F$ is the Hodge filtration induced from ${\cal L}_X$,
  is a polarised variation of Hodge structures $VHS$.

  \smallskip
\n  v) The following decomposition property of $({\Omega}^* {\cal
L}, {\cal W} ) $ (2) into
 intermediate extensions of polarised $VHS$ is satisfied}

 \smallskip
\noindent  \centerline {$(Gr^{\cal W}_{r}{\Omega}^* {\cal L}, F)
\cong \oplus_{K\subset I} {\bf  j}_{!*}^K {\cal L}^K_{r} [ - \vert
K \vert], F[-\vert K \vert])$, for $ r > 0$ }

\smallskip
\centerline {$(Gr^{{\cal W}}_{0}{\Omega}^* {\cal L}, F) \cong 0$}

\smallskip
\centerline {$(Gr^{{\cal W}}_{r}{\Omega}^* {\cal L}, F) \cong
\oplus_{K \subset I} j_{!*}^K {\cal L}^K_{r} [1 - \vert K \vert ],
F) $, for $ r < 0$}

\medskip \n {\bf Theorem 2}. {\em  There exists a rational weight
filtration $({\bf j}_{*}{\cal L}, {\cal W})$ and a
quasi-isomorphism\\ $({\bf j}_{*}{\cal L}, {\cal W}) \otimes_{\Q}
\C \cong ({\Omega}^* {\cal L},{\cal W} )$, such that ${\cal W}_0
\cong {\bf j}_{!*}{\cal L}$.\\}

\medskip \n {\bf Theorem 3}. {\em
For $X$ proper of dim.$n$, there exists a natural $MHS$ on various
cohomology groups with coefficients in ${\cal L}$ underlying a
polarised $VHS$ of weight $m$ on $X-Y$, as follows

\smallskip
\noindent i) The bi-filtered complex

\smallskip
\noindent \centerline {$({\Omega}^*{\cal L}[n], W, F)$ where $W =
{\cal W}[m+n]$  ( that is $W_{i+m+n} = {\cal W}_i $)  for $i \geq
0$ and $0$ for $i < 0 $ }

\smallskip
\noindent  underlies a $MHC$ isomorphic to ${\bf  j}_{*}{\cal
L}[n]$ s.t.  $H^{i}(X-Y, {\cal L}[ n ])$ is of weight $ \geq i + m
+ n$. The filtration $W$ is by perverse sheaves. \\
 Dually $H_c^{i}(X-Y, {\cal L}[ n ])$ is of
weight $ \leq   i + m + n$.}

\smallskip
\noindent {ii)}{\, \em Let $i_Y : Y \to X$ denotes the embedding.
The   quotient complex
 $i_Y^* ({\Omega}^*{\cal L}[n]/ W_{m+n}) $ with the
induced filtrations is a $MHC$ quasi-isomorphic to ${\bf i}_Y^{!}
{\bf j}_{!*}{\cal L}[n+1]$ s.t. $ H^{i}_Y(X, {\bf j}_{!*}{\cal
L}[ n ])$ is of weight $ \geq i + m + n$. }

\smallskip
\noindent {iii)}{\, \em The bi-filtered complex

\smallskip
\noindent \centerline {$(i^{*}_Y {\cal W}_0 {\Omega}^*{\cal L}[n], W,
F)$  where $W = {\cal W}[m+n+1], \, ( W_{i+m+n+1} = {\cal W}_i )$
for $i < 0$}

\smallskip
\noindent  underlies a $MHC$ isomorphic to $i^{*}_Y {\bf j}_{!*}{\cal
L}[n]$ s.t. $ H^{i}(Y, i^{*}_Y {\bf j}_{!*}{\cal L}[ n ])$ is of
weight $ \leq   i + m + n$.}

\smallskip
\n As we see the Intersection cohomology [17] is the fundamental
ingredient which provides the new class of Hodge complexes not
merely defined by complete non singular projective varieties.
 In a different direction, the theory
  has been successfully related to the
theory of differential modules [29].  We give now a more detailed
discussion of the contents.

\smallskip
  \noindent {\it Weight filtration}. The defining  property  of such filtration
  can be understood after a digression
on the sheaves of nearby cycles. Under the additional hypothesis
of the existence of  a local equation $f$ of $Y$, defining  a
morphism  on an open subset $f: U \to D $ to a complex disc, such
that $ Y \cap U = f^{-1}(0)$,
 {\em let ${\cal N} = Log{\cal T}^u$
denotes the logarithm of the unipotent part of the monodromy. The
 filtration  $W({\cal N})$ on  ${\Psi}^u_f$ is defined by
the nilpotent endomorphism  ${\cal N}$ in the abelian category of
 rational (resp.  complex) perverse sheaves.}
The  isomorphism in the abelian category of perverse sheaves [1,2]
$${\bf j}_{*}{\cal L}[ n ]/
 {\bf j}_{!*}{\cal L}[ n ] \; \simeq \; Coker {\cal N}
\colon {\Psi}_f  ^u {(\cal L)} [ n - 1 ] \rightarrow {\Psi}_f^u
{(\cal L)} [ n - 1 ] \leqno{(4)} $$
  \n suggest to start the  weight filtration
with ${\bf j}_{!*}{\cal L}[ n ]$ and continue with $W({\cal N})$
induced on $Coker {\cal N}$, then the main problem is to show,
that {\it the various weights defined locally on $Y$ (for
different local equations) glue together on} ${\bf j}_{*}{\cal L}[
n ]/
 {\bf j}_{!*}{\cal L}[ n ]$.\\
 For each local equation $f$, ${\Psi}_f^u$ can be defined by a
section of Verdier's  monodromic specialization sheaf on the
normal cone to $Y$ in $X$ [11]. For a different equation $g = u f$
with $u$ invertible ${\Psi}_g^u$ and ${\Psi}_f^u$ are isomorphic
but not canonically, so we cannot define a global complex
${\Psi}_Y^u$ on $Y$ and $W({\cal N})$ on it, however we prove that
the isomorphism becomes canonical on $Coker {\cal N}$ hence  the
induced $W({\cal N})$ is globally defined on ${\bf j}_{*}{\cal L}[
n ]/ {\bf j}_{!*}{\cal L}[ n ]$ .\\
 Here we solve the problem in the case of
$NCD$ by a different approach. We define first the complex weight filtration and then
prove it is rationnally defined. The detailed study of the complex
$\Omega^*{\cal L}[n]$, necessary for Hodge theory, leads to such
approach of the rational weight.

\smallskip \n {\bf Theorem 4}. {\em For each local equation $f$ of
$Y$ defined on an open set $U$ of $X$, ${\cal N}$ acts on
${\Psi}^u_f {\cal L}[n-1]$; then
  the induced filtration by ${\cal W}({\cal N})$ on  $Cokernel$ of
  ${\cal N}$ ( the right term of
(4)) coincides with the induced filtration by ${\cal W}$ on  $
 {\bf j}_{!*}{\cal L}[n] / {\bf j}_{*}{\cal L}[n]$ (the
 left term)}.

\smallskip
\n This problem has been solved successfully for any divisor using
the deformation to the normal cone  and will appear later in a
joint work with L\^e D.T. and Migliorini L.

 \smallskip
\n {\em Let $p: C_YX - Y \to Y$ denotes the projection
 of the punctured normal cone onto $Y$. There exists
 a canonical filtration ${\cal W}$  on ${\bf
j}_*{\cal K}$ defined as  ${\cal W}_0 = { j}_{!*}{\cal K}$ and for
$i > 0$ \\
${\cal W}_i ({\bf j}_{*}{\cal K}/
 {\bf j}_{!*}{\cal K} )\; \simeq \;{\cal W}_{i-1}{^p{\cal
H}^0}({\bf p}_*(\Psi_Y{\cal K}))\colon= Ker [{^p{\cal H}^0}({\bf
p}_* (\Psi_Y{\cal K}) \to {^p{\cal H}^0}({\bf p}_* (\Psi_Y{\cal
K}/{\cal W}_i({\cal N}{\vert
\Psi_Y{\cal K})]} $\\
s.t. for each local equation $f = 0$ of $Y$  defined on an open
set $U_f$, we have on $Y_f = U_f \cap Y$\\
${\cal W}_i ({\bf j}_{*}{\cal K}/
 {\bf j}_{!*}{\cal K} )_{\vert Y_f}\; \simeq \;
{\cal W}_{i-1}( Coker {\cal N}\vert ({\Psi}^u_f  {\cal K}[-1])$.

\smallskip \n Dually on $i^*{\bf j}_{!*}{\cal K} [-1]\simeq
{^p{\cal H}^{-1}}({\bf p}_* (\Psi_Y{\cal K})$ the weight is
defined as\\
 $({\cal W}_i \, i^*{\bf j}_{!*}{\cal K}) [-1]
 \simeq {^p{\cal H}^{-1}}({\bf p}_* {\cal W}_{i+1}({\cal
N}{\vert \Psi_Y{\cal K})})$ }
===========
\smallskip
\noindent {\it Purity and decomposition (local results)}. Working
with the complex local system ${\cal L}^c$ we need to exhibit a
bifiltered complex $({\Omega}^*{\cal L}, {\cal W}, F)$ underlying
a mixed Hodge complex ($MHC$). Since the local system itself
results in the geometric case from singularities of morphisms, it
is natural to search for a combinatorial logarithmic complex for
such bi-filtered complex. Its construction is suggested by an
algebraic formula of the Intersection complex given by Kashiwara
and Kawai in [26].
 We  explain now
the main  basic local results.\\
1- If we consider a point $y \in Y_M^*$, a variation of Hodge
structures $VHS$ on ${\cal L} $ of weight $m$ defines  a nilpotent
orbit $L$ with a set of nilpotent endomorphisms $N_i, i \in M $.
The nilpotent orbit theorem [4], [24] states that the $VHS$
degenerates along $Y_M^*$ into a variation of $ MHS $ with weight
filtration $W^M = W(\Sigma_{i\in M} N_i)$ shifted by $m$.\\
However the difficulty in the construction of the weight is to
understand what happens  at the intersection of $Y_M$ and $Y_K$
for two subsets $M$ and $ K$ of $I$. This difficulty  couldn't be
explained for the $\Q-$ structure until the discovery of perverse
sheaves.
 In order to prove the decomposition of $Gr^{\cal
W}{\Omega}^*{\cal L}$ into intermediate direct image of various
local systems on the components of $Y = \cup_{i\in I} Y_i $, we
introduce in (3.3) for $K \subset M$, the complex $C^K_r L$
constructed  out of the nilpotent orbit $(L, N_i)_{i\in K}$
defined by ${\cal L}$ at a point $z \in Y^*_K $ and $W^K = W(
\Sigma_{i\in K} N_i)$. We  prove that $C^K_r L$ has a unique non
vanishing cohomology isomorphic to $Gr^{W^K}_{r - \vert K \vert}
(L/\Sigma_{i\in K} N_i L)$,  the fiber at $z$  of the rationally
defined local system ${\cal L}^K_r$
previously introduced  for $r>0$.\\
Now $Y_M $ is a subset of $Y_K$ and we are interested in the fiber
of $j_{!*}^K {\cal L}^K_r$ at $y \in Y_M^*$, so we introduce the
complex $C^{KM}_r L$ in (3.3), which is quasi-isomorphic to this
fiber and will appear as a  component of the decomposition (3.8)
of the graded part of the weight filtration.\\
These main results in the open case form the content of  the
second and third sections.\\
 2- The second major  technical result (5.3) in the case of nearby
 co-cycles, shows that if
$f$ is a local equation of $Y$, $C^K_r L$ is  isomorphic to the
primitive part  of   the fiber $\nu^K_r L$ at $y$ of a local
system on $Y_K^*$ component of the decomposition of
 $Gr_r^{W({\cal N})}{\psi}_f^u{\cal L}$, which proves that the
 complex weight filtration is compatible with our description of
 the rational weight.\\
3- The key result that enables us to give most of the proofs is
the existence, for a nilpotent orbit ($L, N_i, i \in K$), of a
natural decomposition  ( 2.4 )\\
 \centerline { $Gr^{W^K}_r  L  = \oplus Gr^{W(N_{i_n})}_{m_{i_n}} \cdots
Gr^{W(N_{i_1})}_{m_{i_1}} L \, : \, \Sigma_{i_j\in K} m_{i_j} = r
$}

\n As a consequence the   local systems on $Y^*_K$ defined for $r
> 0$ by $C^K_rL$ decompose into a direct sum of elementary
components\\
 \centerline {$ \oplus Gr^{W(N_{i_k})}_{m_{i_k}} \cdots Gr^{W(N_{i_1})}_{m_{i_1}} (L /
(N_{i_1} L + \cdots + N_{i_k} L)) \, : \, \Sigma_{i_j \in K}
m_{i_j} = r + \vert K \vert, m_{i_j} \geq 2 $}

\n the corresponding elementary complexes are introduced in (3.5)
and are key
ingredients in the proof.\\
The above results explain the subtle relation between  the
filtration $\cal {W} (\cal {N})$ on $\Psi^u_f ({\cal L})$ which is
hard to compute and the various local monodromy at points of $Y$.

\smallskip
\noindent We use the terminology of  perverse sheaves but we leave
to the reader
 the choice in the shift in degrees, generally  by $-n$. \\
Finally the local definition of the weight is in $(3.2)$, purity
in (3.6) and decomposition in (3.8). The global  definition of the
weight is in $(4.1)$ and the decomposition in (4.3).
 The weight of the nilpotent action on $\Psi^u_f$ is in $(5.1)$
and the comparison
 in $(5.3)$. Finally we suggest strongly to
the reader to follow the proofs on an example, sometimes on the
surface case as in the example in $(3.2)$; this example will be again useful for
    $\Psi_f {\cal L}$ in $(5.3)$. For $X$ a line and $Y = 0$ a point,
 the fiber at $0$
 of ${\bf j}_{*}{\cal L} $ is a complex $ L \stackrel{N}{\rightarrow} L$
  where $(L, N)$ is a nilpotent orbit  of weight $m$ and the weight ${\cal W}[m]$
   on the complex  is $({\cal W}[m])_{r+m} = {\cal W}_r $ defined by
   the sub-complex $ (W_{r+1}L \stackrel{N}
   {\rightarrow} W_{r-1}L)$.
\section{ Local invariants of  ${\cal L}$}

We  need a precise description in terms of the local invariants of
the local system $ {\cal L}$. We recall some preliminaries on
${\bf j}_{*}{\cal L} $ and  we give a basic local decomposition of
weight filtrations defined by the local monodromy on nilpotent
orbits.
\subsection{ Preliminaries }

 \be { Local and global description of ${\bf
j}_{*}{\cal L} $.}{1.II.1}

\n In the neighbourhood of a point $y$ in $Y$, we can suppose $X
\simeq D^{n+k}$ and $X^* = X-Y \simeq {(D^*)}^n \times D^k$ where
$D$ is a complex disc, denoted with a star when the origin is
deleted. The fundamental group $\Pi_1 (X^*)$ is a free abelian
group generated by $n$ elements representing  classes of closed
paths around the origin, one for each $D^* $ in the various axis
with one dimensional  coordinate $z_i$  ( the hypersurface $Y_i$
is defined locally by the equation $z_i=0$ ). Then the local
system $\cal L$ corresponds to a representation of $ {\Pi}_1(X^*)$
in a vector space  $L$ defined by the action of commuting automorphisms
$T_i$ for $i \in [1,n ]$ indexed by the local components $Y_i$ of
$Y$ and called monodromy action around $Y_i$. The automorphisms
$T_i$ decomposes as a product of commuting automorphisms,
semi-simple
and unipotent $T_i = T_i^s T_i^u$.\\
Classically $L$ is viewed as the fiber of $\cal L$ at the
reference point for the fundamental group $\Pi_1 (X^*)$, however
since we will need  to extend the Hodge filtration on Deligne's
extended bundle, it is important to view  $L$ as the vector space
of multiform sections of $\cal L$ (that is the sections of the
inverse of  $\cal L$  on a universal covering of $X^*$ ).

\smallskip
 \noindent Given a $\Q-$local system, locally unipotent
  along  $Y$ (to simplify the exposition) we consider ${\cal L}_{X^*} \colon=
 {\cal L} \otimes_{\Q} {\cal O}_{X^*}$ and its Deligne's bundle extension ${\cal L}_X$
 which has a nice  description as the subsheaf of ${\bf j}_* {\cal L}_{X*}$
generated  locally at a point $y$ in $Y$ by sections associated to
multiform sections of ${\cal L}$ as follows. \\
The logarithm of the unipotent monodromy, $ N_i \colon= Log T_i^u
= -{\Sigma}_{k \geq 1} (1/ k) {(I-T^u_i)}^k $ is  defined as the
sum of nilpotent endomorphisms $(I - T^u_i)$ so that this sum is
finite.\\
 A multiform section  $v$
 corresponds to a germ $\mathaccent"7E{v} \in
{\bf j}_* {\cal L}_{X*}$ with an explicit description of the
action of the connection by the formulas
$$ \mathaccent"7E{\hbox{v}} (z)  = (exp (-\frac{1}{2 i \pi}
\Sigma_{j \in J} ( log
z_{j}) N_{j})). v, \quad \nabla \mathaccent"7E {v} = -\frac{1}{2 i \pi}\Sigma_{j \in J}
 \widetilde { N_{j}.v } \otimes \frac
{dz_j}{z_j} \leqno{(5)}$$
 \noindent a basis of $L$ is sent on a basis of
${\cal L}_{X,y}$.\\
 The residue of the
connection $\nabla $ along each  $Y_j$ defines an endomorphism $
{\cal N }_j $ on the restriction ${\cal L}_{Y_j} $ of ${\cal
L}_X$.

\smallskip
\n  The fiber at the origin of the complex ${\Omega}^*_X (Log Y)
\otimes {\cal L}_X $  is quasi-isomorphic to a Koszul complex as
follows. We associate to $(L, N_i), i \in [1,n ]$   a strict
simplicial vector space such that for all sequences $(i.) = (i_1 <
\cdots < i_p)$

\smallskip
\centerline {$L(i.) = L \qquad , \qquad N_{i_j} \colon L(i. - i_j)
\rightarrow L(i.)$}

\smallskip
\n {\bf Definition.} {\em The  simple complex defined by the
simplicial vector space above is the Koszul complex (or the
exterior algebra) defined by $(L, N_i)$ and denoted by $\Omega (L,
N.)$. A general notation is  $ s(L(J), N.)_{J \subset [1,n ]} $
 where  $J$ is identified with the strictly increasing
sequence of its elements and where $ L (J) = L $}.

\smallskip
\n  It is quasi-isomorphic to the Koszul complex $\Omega (L, Id -
T.)$ defined by $(L, Id - T_i),  i \in [1,n ]$. This local setting
compares to the global case via Grothendieck and Deligne DeRham
cohomology results.

\smallskip
\n {\bf Lemma}. {\em  For $M \subset I$ and  $y \in Y^{*}_M$, the
above correspondence $v \mapsto \tilde {v} $, from $L$ to  ${\cal
L }_{X,y}$,  extended from $L(i_1, \ldots, i_j)$ to
$({\Omega}^*_X (Log Y) \otimes {\cal L}_X)_y $ by $v \mapsto
\tilde {v} \frac {dz_{i_1}}{z_{i_1}}\wedge \ldots \wedge  \frac
{dz_{i_j}}{z_{i_j}}$, induces quasi-isomorphisms
$$({\Omega}^*_X (Log Y) \otimes
{\cal L}_X)_y \cong
 \Omega (L,  N_j, j \in M )\cong s(L(J), N.)_{J \subset
[1,n ]} \leqno{(6)}$$
  hence $ ({\bf j}_* {\cal L}^c)_y \cong
 \Omega (L,  N_j, j \in M) $.}

\smallskip
\n  This description of  $({\bf j}_* {\cal L })_y$ is the model
for the description of the next various perverse sheaves.
 \be{The intermediate extension ${\bf j}_{!*}{\cal L} $}{1.II.2} \
   Let $ N_J = \Pi_{j \in J}  N_j$ denotes a composition
   of endomorphisms of $L$, we consider  the strict  simplicial  sub-complex
of the DeRham logarithmic complex  defined by $Im N_J$ in $L(J) =
L$.

\smallskip
\n {\bf Definition}. {\em The  simple complex defined by the above
simplicial sub-vector space  is  denoted by}
$$IC (L)\colon = s(N_J L, N.)_{J \subset M}, \ N_J L
\colon= N_{j_1} N_{j_2}\ldots N_{j_p}L , j_i \in J$$ Locally the
germ of the intermediate extension ${\bf j}_{!*} {\cal L}$ of
${\cal L}$   at a point $y \in Y^*_M$  is quasi-isomorphic to the
above complex [24 (3)]
$$ {\bf j}_{!*} ({\cal L})_y \simeq IC
(L) \simeq s(N_J L,  N.)_{J \subset M} $$ The corresponding global
DeRham description is given as a sub-complex $IC (X, {\cal L})$ of
 ${\Omega}^*_X (Log Y) \otimes {\cal L}_X$. In terms of a set of
  coordinates $z_i, i\in M, $
defining $Y_M$
 in a neighbourhood of $y \in Y^*_M$,
 $IC^j (X, {\cal L})$  is the subanalytic complex
 of ${\Omega}^*_X (Log Y) \otimes {\cal L}_X$ with fiber at $y$
  generated, as an ${\Omega}^*_{X,y}$
sub-module, by the sections $\tilde {v} \wedge_{j \in J} \frac
{dz_j}{z_j}$ for $v \in N_J L$. This formula is independent of the
choice of coordinates, since if we choose a different coordinate
$z'_i = f z_i$ instead of $z_i$, with $f$ invertible holomorphic
at $y$, the difference $ \frac {dz'_i}{z'_i} - \frac {dz_i}{z_i} =
\frac {df}{f}$ is holomorphic at $y$, hence the difference of the
sections $\tilde {v} \wedge_{j \in J} \frac {dz'_j}{z'_j} - \tilde
{v} \wedge_{j \in J} \frac {dz_j}{z_j}$  is still a section of the
sub-complex  $IC (X, {\cal L})$; moreover the restriction of the
section is still defined in the sub-complex near $y$, since $N_J L
\subset N_{J-i} L$ for all $i \in J$.

\smallskip
\n {\bf Lemma.} We have a quasi-isomorphism   ${\bf j}_{!*}{\cal
L}^c \cong IC (X, {\cal
 L})$.
 \be {Hodge filtration and Nilpotent orbits  }{1.II.3}
 \n {\it  Variation
of Hodge structures (VHS).} Consider the flat bundle $({\cal L}_X,
\nabla) $ in the previous hypothesis and suppose now that ${\cal
L}_{X^*}$ underlies a $VHS$ that is a polarised filtration by
sub-bundles $F$ of weight $m$ satisfying Griffith's conditions
[19].

\smallskip
\n The {\it nilpotent and the $SL_2$ orbit  } theorems [19], [4],
[24], [25] show that $F$ extends to a filtration by sub-bundles
$F$ of ${\cal L}_X$ such that the restrictions to open
intersections $Y^*_M $ of components of $Y$ underly {\it locally}
a variation of mixed  Hodge structures $VMHS$ where the weight
filtration is defined by the nilpotent endomorphism ${\cal N}_M$,
residue of the connection, (there is no flat bundle defined
globally on $Y^*_M $, if $z_{i_1}, \ldots,z_{i_n}$
for $i_j \in M$ are local equations at $y \in Y^*_M$,
then $\Psi_{z_{i_n}} \circ \ldots \circ \Psi_{z_{i_1}}{\cal L} $ is
the underlying local system nar $y$).

\smallskip
\n {\it Local version.} Near a point $y \in Y^*_M$ with $\vert
M\vert = n$ a neighbourhood of $y$ in the fiber of the normal
bundle looks like a disc $D^n$ and the above hypothesis reduces to

\smallskip
\noindent {\it Local Hypothesis :  Nilpotent orbits [4].} Let
$$ (L, N_i, F, P, m, i\in M = [1,n] ) \leqno{(7)} $$
be defined by the $VHS$, that is a $\Q-$vector space $L$ with
endomorphisms $N_i$ viewed as defined by the multiform horizontal
(zero) sections of the connection on $ (D^*)^n$ (hence sections on
the inverse image on the universal covering), a Hodge structure
$F$ on $L^{\C} = L \otimes_{\Q} \C $ viewed as the fiber of the
vector bundle ${\cal L}_X $ at $y$ (here $y = 0 $), a natural
integer $m$ the weight and the polarization $P $.\\
The main theorem in [4] states that   for all
$N = \Sigma_{i \in M} {\lambda}_i N_i$
with ${\lambda}_i > 0$ in ${\hbox {\bf R}}$ the
 filtration $W(N)$ ( with center
$0$ ) is independent of $N$ when  ${\lambda}_i$ vary and $
W(N)[m]$ is the weight filtration of a  graded polarised $MHS$
 called the  limit $MHS$  of weight $m$ ( $L, F,W(N)[m]) $.

\smallskip
\noindent {\em Remark: $ W(N)[m]$ is $ W(N)$ with indices shifted by $m$ to the
right: $ (W(N)[m])_r\colon= W_{r-m}(N)$, the convention being a
shift to left for a decreasing
filtration and to right for an increasing filtration.}

\smallskip
\noindent It is important to notice that the  orbits depend on the
point $z$ near $y$ considered, in particular $F_z \neq F_y $. In
this case when we restrict the orbit to $J \subset M $, we should
write\\
\centerline {   $ (L, N_i, F (J), P, m, i\in J \subset M ) $}

\n We write $W^J $ for $ W(N_J) $ where $ N_J = \Sigma_{i \in J}
 N_i$. We will need the following result [4 p 505]:\\
 Let $I, J \subset M$
 then $W^{I\cup J} $ is the
 weight filtration of $N_J$ relative to $W(N_I)$\\
\centerline { $ \forall j,i \geq 0, N_J^i: Gr^{W^{I\cup J}}_{i+j}
Gr^{W^I}_j {\buildrel \sim \over \rightarrow } Gr^{W^{I\cup
J}}_{j-i} Gr^{W^I}_j.$}
\subsection{Properties of the relative weight filtrations}
Given a nilpotent orbit  we may consider various filtrations $W^J
= W(\Sigma_{i\in J}N_i)$ for various $J \subset M$. They are
centered at $0$ ( that is  we suppose here the weight of the
nilpotent orbit equal to zero, otherwise  the true weight of the
$MHS$ is defined up to a shift), preserved by $N_i$ for $i\in M$
and shifted by $-2$ for $i \in J$: $N_iW^J_r \subset
N_iW^J_{r-2}$. We need to know more about the action of $N_i$
which is compatible with $W(N_j)$. The starting point of this
study is the definition of the relative weight filtration by
Deligne [10] \\ {\em Let $(L, W)$ be  a finite dimensional vector
space $L$ endowed with an increasing filtration $W$ and $N$ a
nilpotent endomorphism compatible with $W$. There may exists at
most
a unique filtration $M = M(N,W)$ satisfying \\
1) $ N: M_j \simeq M_{j-2}$\\
2) $ N^j: Gr_{k +j}^{M} Gr_k^{W}L \simeq
 Gr_{k-j}^{M} Gr_k^{W}L,$}\\
A main result in [4] shows \\
{\em The filtrations $W^J$
defined by a polarised nilpotent orbit satisfy\\
1) For a subset $J \subset [1,n]$, $ \forall  j \in J,  \forall \lambda_j >
0 $, $N_J =
 \sum_{j \in J} \lambda_j N_j$, the filtration
$W^J = W(N_J)$ is independent of $\lambda_j > 0 $ \\
2) For subsets  $J$ and $J'$ in $[1, n]$,
$A = J \cup J'$ we have for all $j \in \N, k \in \Z$ :\\
 $ N_J^j: Gr_{k +j}^{W^A} Gr_k^{W^{J'}}
L \simeq
 Gr_{k-j}^{W^A} Gr_k^{W^{J'}}L$,\\
 that is $W^A$ is the relative weight filtration of $N_J$ acting
 on $ (L, W^{J'})$}.\\
 We remark  also that  for all $ J \subseteq B \subseteq A$,
$N_B$ and $N_J$ induce the same morphism on
 $Gr_k^{W^{J'}}L$.\\
 Finally we need the following result of Kashiwara
 ( [25, thm 3.2.9, p 1002])\\
  {\em Let $(L,N, W) $
consists of a vector space endowed with an increasing filtration
$W$ preserved by a nilpotent  endomorphism $N$ on $L$ and suppose
that the relative filtration $ M = M(N,W)$ exists, then there
exists a canonical decomposition:}  $$Gr^M_l L = \oplus_k Gr^M_l
Gr^W_k L
$$ Precisely, Kashiwara exhibits a splitting of  the exact sequence: \\
$0 \rightarrow W_{k-1} Gr^M_l L \rightarrow W_k Gr^M_l L
\rightarrow Gr^W_k Gr^M_l L \rightarrow 0$. \\
 by constructing a natural section of $ Gr^W_k
Gr^M_l L$ into $ W_k Gr^M_l$. We will need later more precise
relations between these filtrations that we discuss now.
  \be {Key lemma}{2.4}{\em
(Decomposition of the relative weight filtrations) : Let $ (L,
N_i, i \in [1,n], F)$ be a polarised nilpotent orbit and for $A
\subseteq [1,n]$ let $W^A \colon= W(\sum_{i \in A}
 N_i)$ (all weights centered at $0$ ), then
:\\
i) For  all $ i \in A$, the filtration $W^A$ induces a
 trivial filtration on $Gr_k^{W^i} Gr_{k'}^{W^{(A-i)}}L $ of weight
 $k+k'$ \\
 ii) For  $ A = \{ i_1,\ldots,
i_j\}\subset [1,n]$, of length  $\vert  A \vert  = j $ we have a
natural decomposition
$$ Gr_r^{W^A} L \simeq \oplus_{m_{i.} \in X^A_r} Gr^{W^{i_j}}_{m_{i_j}}
\cdots Gr^{W^{i_1}}_{m_{i_1}} L \;{\hbox {where}}\; X^A_r =
\{m_{i.} \in {\Z}^j: \Sigma_{i_l \in A } m_{i_l} = r\} $$ more
precisely
$$ Gr_r^{W^A} (\cap_{i_s \in A} W^{i_s}_{a_{i_s}} L )
\simeq \oplus_{\{m_{i.}\in X^A_r, m_{i_s} \leq a_{i_s}\}}
Gr^{W^{i_j}}_{m_{i_j}} \cdots Gr^{W^{i_1}}_{m_{i_1}} L \leqno{(8)}
$$
 \noindent  iii) Let $A = B \cup C $,  $ N_i' $ denotes the
restriction of $N_i$ to $Gr^{W^C}_c$ and $ N_B' = \Sigma_{i\in
B}N'_i$, then  $W_b^B $ induces $W_b( N_B')$ on $Gr^{W^C}_c$, that
is
$$ Gr_{b+c}^{W^A} Gr_c^{W^C} L \simeq Gr_b^{W^B} Gr_c^{W^C} L \simeq
Gr_{b+c}^{W^A} Gr_b^{W^B} L \simeq  Gr_c^{W^C} Gr_b^{W^B} L $$
\noindent iv) The repeated graded objects in i) do not depend on
the order of the elements in $A$.}

\smallskip

\noindent Remark: This result give relations between  various
weight filtrations in terms of the elementary ones $W^i \colon=
W(N_i)$ and will be extremely useful in the study later of the
properties of the  weight filtration on the mixed Hodge complex.

\smallskip
\n {\it Proof }. To stress the properties of commutativity of the
graded operation for the filtrations, we prove first

\smallskip
\noindent {\em Sublemma:  For all subsets $[1,n]\supset A \supset
\{B,C\},$ the isomorphism of Zassenhaus $Gr^{W^B}_b Gr^{W^C}_c L
\simeq Gr^{W^C}_c Gr^{W^B}_b L $ is an isomorphism of $MHS$ with
weight filtration (up to a shift) $W = W^A$ and Hodge filtration
$F =  F_A $, hence compatible with the third filtration $W^A$ or
$F_A$.}

\smallskip
\n Proof of the sub-lemma: Recall that both spaces $Gr^{W^B}_b
Gr^{W^C}_c $ and $ Gr^{W^C}_c Gr^{W^B}_b$ are isomorphic to $W^B_b
\cap W^C_c $ modulo $ W^C_c \cap W^B_{b-1} + W^B_b \cap W^C_{c-1}
$. In this isomorphism a third filtration like $F_A$ (resp. $ W^A
$) is induced on one side by $ F'_k = (F_A^k \cap W^C_c) +
W^C_{c-1}$ (resp. $W'_k = (W_k^A \cap W^C_c) +  W^C_{c-1}$) and on
the second side by  $ F_k'' = (F_A^k \cap W^B_b) +  W^B_{b-1}$
(resp. $W_k'' =
 W_k^A \cap W^B_b) +  W^B_{b-1}$). We
introduce the third  filtration $ F_k''' = F_A^k \cap W^B_b \cap
W^C_c $ (resp. $W_k''' = W_k^A \cap W^B_b \cap W^C_c $) and we
notice that all these spaces are in the category of $MHS$, hence
the isomorphism of Zassenhaus is strict and compatible with the
third filtrations induced by $F_A$ ( resp $W^A$ ).

\smallskip
\n Proof of the key lemma i). Let $ A \subset [1, n]$ and $i\in
A$, then $W^A $ exists on $L$ and coincides with  the relative
weight filtration for $N_i$ with respect to $W^{(A-i) }$ by a
result of Cattani and Kaplan. Then we have by Kashiwara's result $
Gr_l^{W^A}L \simeq \oplus_k Gr_l^{W^A} Gr_k^{W^{(A-i)}}L  $. Let
us attach to each point $(k,l)$ in the plane the space $
Gr_l^{W^A} Gr_k^{W^{(A-i)}}L  $ and let $M_j = \oplus_l Gr_l^{W^A}
Gr_{l-j}^{W^{(A-i)}}L $ be the direct sum along indices in the
plane $(k,l)$ on a parallel to the diagonal ($l = k + j$). Then we
have for $j>0$
$$ (N_i)^j: Gr_{k +j}^{W^A} Gr_k^{W^{(A-i)}}L \simeq
 Gr_{k-j}^{W^A} Gr_k^{W^{(A-i)}}L, \; (N_i)^j: M_j \simeq M_{-j}.$$
This property leads us to introduce the space  $V = \oplus_l
Gr_l^{W^A } L \simeq \oplus_{l,k}  Gr_l^{W^A} Gr_k^{W^{(A-i)}}L $,
 then $N_i$ on $ L$ extends to a nilpotent endomorphism on $V$,
  $ N_i: V \rightarrow V $ inducing $ N_i:Gr_l^{W^A}L
\rightarrow Gr_{l-2}^{W^A }L $ on each $l-$component of $V$. We
consider on $V$ two increasing filtrations  $W'_s \colon=
\oplus_{l-k \leq s} Gr_l^{ W^A }Gr_k^{W^{(A-i)}}L$ and
  $ W_s'' \colon= \oplus_l W_s^i Gr_l^{W^A}L $. Then
$ N_i$ shift these filtrations by $-2$. In fact $N_i: W'_s
\rightarrow W'_{s-2}$ sends $Gr^{W^A}_lGr_k^{W^{(A-i)}}L$ to
$Gr^{W^A}_{l-2}Gr_k^{W^{(A-i)}}L$ and $ (N_i)^j $ induces an
isomorphism $ Gr^{W'}_j V \simeq Gr^{W'}_{-j} V$. As well we have
an isomorphism $ Gr^{W''}_j V \simeq Gr^{W''}_{-j} V$, since
$(N_i)^j : (Gr^{W^i}_j L, W^A, F_A) \simeq  (Gr^{W^i}_{-j} L, W^A,
F_A) $ is an isomorphism of $MHS$ up to a shift in indices, hence
strict on $ W^A $ and $ F_A$ and induces an isomorphism  $
Gr^{W^i}_j Gr_l^{W^A} \simeq Gr^{W^i}_{-j}Gr_{l-2j}^{W^A} $.
  Then these two filtrations $W'_s $ and $W_s''$ are equal by uniqueness
of the weight filtration of  $ N_i$ on $V$, that is
$$ W_s'' = \oplus_{k,l} W^i_s Gr_l^{W^A } Gr_k^{W^{(A-i)}}L = W'_s =
\oplus_{l-k \leq s} Gr_l^{W^A} Gr_k^{W^{(A-i)}}L $$ that is $
W^i_s Gr_l^{W^A } Gr_k^{W^{(A-i)}}L =  Gr_l^{W^A }
Gr_k^{W^{(A-i)}}L $ if $l-k \leq s$ and $ W^i_s Gr_l^{W^A }
Gr_k^{W^{(A-i)}}L = 0 $ otherwise, hence
$$ Gr_l^{W^A} Gr^{W^i}_j L \simeq  Gr_l^{W^A} Gr_j^{W^i}
Gr_{l-j}^{W^{A-i}}L, {\hbox {and for all}} \; l \neq k + k', \;
Gr_l^{W^A} Gr_k^{W^i} Gr_{k'}^{W^{(A-i)}}L \simeq 0. $$
 which ends the proof of (i).

 \n ii) Since $W^A$ induces a trivial filtration on $Gr_k^{W^i}
Gr_{k'}^{W^{(A-i)}}L $ of weight $k+k'$ we have
$$ Gr_l^{W^A} L \simeq \oplus_k Gr_l^{W^A} Gr_k^{W^{A-i}} L  \simeq
\oplus_k Gr_l^{W^A} Gr_k^{W^{A-i}} Gr^{W^i}_{l-k} L \simeq
\oplus_k Gr_k^{W^{(A-i)}} Gr_{l-k}^{W^i} L.$$ Now if we suppose by
induction on length of $A $, the decomposition true for $ A-i$, we
deduce easily
 the decomposition for $ A$ from the above result.

\noindent iii) We restate here the property of the relative
 monodromy for $W^A$ with respect to $W^C$ and we apply ii).

\noindent iv)In the proof above we can start with any $i$ in $A$,
hence the decomposition is symmetric in elements in $A$. It
follows that the graded objects of the filtrations $W^i, W^r,
W^{\{i,r,j\}}$ commute and since $W^j$ can be expressed using
these filtrations, we deduce that $W^i, W^r, W^j$ also commute,
for example: $ Gr_{a+b+c}^{W^{\{i,r,j\}}} Gr_{a+b}^{W^{\{i,j\}}}
Gr_a^{W^i} \simeq Gr_c^{W^r}Gr_b^{W^j} Gr_a^{W^i}$ is symmetric in
$i,j,r$.
\section{The weight filtration and main  theorems in the local case}
To describe the weight filtration, we introduce a category $S(I)=
S$ attached to a set $I$ already used by Kashiwara and Kawai [26]
for the intersection complex. We start with a local study, that is
to say with the hypothesis of a polarised nilpotent orbit and we
describe the weight filtration $W$ on a combinatorial complex
quasi-isomorphic to the DeRham complex $ \Omega (L,N.)$. The
features of the purity theory will appear relatively quickly.
First we ask the reader to take some time to get acquainted with
the new category $S(I)$ whose objects are indices for the
combinatorial complex. The weights zero or $-1$ describe a complex
$ IC(L)$ quasi- isomorphic to the fiber of the intermediate
extension of ${\cal L}$ and for the other weights we need to
introduce the complexes $C^{KM}_r L$ for $ K \subset M \subset I $
( 3.3)
 which describe the purity theory (3.6) and the
geometry of the decomposition theorem (4.3). A basic technique in the
proof is the decomposition into elementary complexes (3.5),
reflecting
 relations between the weight filtrations of the various
$N_i$.
\subsection{Construction of the weight filtration}
\be {Complexes with indices in the category $S(I)$.} {3.1} \ The
techniques are similar to the simplicial techniques in Deligne's
paper. Here the singularities may come from the coefficients as
well as the $NCD$ in $X$. We introduce a category $S(I)$
attached to a set $ I$, whose objects consist of sequences of
increasing subsets of $I$ of the following form:
$$ (s_.) = (I = s_1 \supsetneqq s_2 \ldots \supsetneqq s_p
\neq \emptyset), \; (p > 0)$$ Subtracting a subset $s_i$ from a
sequence $s_.$ defines a morphism $\delta_i (s.) : (s. -
s_i)\rightarrow s.$ and more generally $Hom ( s.', s.)$ is equal
to one element iff $(s.') \leq (s.)$ that is $(s.')$ is obtained
from $(s.)$ by deleting some subsets. We write $s. \in S(I)$ and
define its degree or
 length $\vert  s.\vert $ as the number of subsets $s_i$
in $(s_.)$.

\medskip
{\em Correspondence with an open simplex. } If $I = \{1, \ldots, n
\}$ is finite, $S(I)$ can be realised as a barycentric subdivision
of the open simplex $\Delta_{n-1}$ of dimension $ n - 1 $. A
subset $K$ corresponds to the barycenter of the vertices in $K$
and a sequence of subsets to an oriented simplex defined  by the
vertices
associated to the subsets.\\
For example, for $I = \{1,2\}$, $S(I)$ consists of the barycenter
$\{3/2 \}$ of $]1,2[$ defined by $ \{1,2\}$, and the open
simplices $ ]1,3/2[$, $]3/2, 2[$ defined resp. by the sequences
$\{1,2\}\supset  \{1\}$ and $\{1,2\}\supset  \{2\}$.\\
 Since all
sequences contain $I$, all corresponding simplices must have the
barycenter defined by $I$ as vertex, that is a sub-simplex
contained in the open simplex $\Delta^*_{n-1} = \Delta_{n-1} -
\partial \Delta_{n-1}$. In this way we define an incidence relation
$\varepsilon(s.,s.')$ between two adjacent sequences equal to $+1$
or $-1$ according to orientation. Incidence relations
$\varepsilon( \Delta_{n-1},s.)$ are defined as well between $
\Delta_{n-1}$ and the simplices corresponding to maximal
sequences.

\smallskip {\it Combinatorial  objects of an abelian category with
indices in $S(I)$ }, that is  functors, are thus
defined, as well as complexes of such objects.\\
 We need essentially the following construction. An algebraic  or analytic
 variety  over a fixed variety $X$ with indices in $S(I)$ denoted by $\Pi$ is a
covariant functor defined by $\Pi( s.): X_{s.} \to X $ and
morphisms $\Pi (s.' \leq s.): X_{s.'} \rightarrow X_{s.}$ over $X$ for
$s.\in S$. An abelian sheaf over $\Pi$ (resp. complex of abelian
sheaves ) ${\cal F}$  is a contravariant functor of abelian
sheaves (resp. complex of abelian sheaves) ${\cal F}_{s.}$ over
$X_{s.}$ (with functorial morphisms $\varphi(s.' \leq s.):
\Pi^* {\cal F}_{s.}\rightarrow {\cal F}_{s.'}$ for
$(s.') \leq (s.)$).

\n The direct image of an abelian sheaf over $\Pi$ (resp. complex
of sheaves ) denoted $\Pi_*{\cal F}$ or preferably  $s({\cal F}_{s.})_{s.\in
S}$ is the simple complex (resp.simple complex associated to a
double complex)  on $X$:
$$s({\cal F}_{s.})_{s.\in
S} \colon= \Pi_*{\cal F} \colon= \oplus_{s. \in S}(\Pi_*{\cal F}_{s.}) [\vert
s.\vert  - \vert  I \vert   ], \; d = \Sigma_{i\in [1,\vert
s.\vert ]} (-1)^{\varepsilon(s.,\delta_i (s.))} \varphi(\delta_i
(s.) \leq s.)$$
 Example. The variety $X$ defines the constant variety
$X_{s.}= X$. The constant sheaf $\Z$ lifts to a sheaf on $X_{s.}$
such that the ''diagonal morphism'' : $\Z_X \rightarrow
\oplus_{\vert (s.)\vert  = \vert  I \vert } \Z_{X_{s.}}$ ( that is
: $n \in \Z \rightarrow (\ldots, \varepsilon( \Delta_{n-1},s.)
n_{s.}, \ldots) \in \oplus_{\vert (s.)\vert = \vert  I \vert } \Z
$ defines a quasi-isomorphism $\Z_X \cong \Pi_*\Pi^* (\Z_X)$. This
is true since $S(I)$ is isomorphic to the category defined by the
barycentric subdivision of an open simplex of dimension  $\vert  I
\vert  - 1$.
  \be { Local definition of the weight filtration}{3.2}
\noindent Our hypothesis here consists again of the polarised
nilpotent orbit $ (L, (N_i)_{i\in M}, F ,m ) $   of weight $m$ and
the corresponding filtrations $ (W^J)_{J\subset M} $ where $W^J =
W (\sum_{i \in  J} N_i)$.

\noindent We will use the category $S(M)$ attached to $M $ whose objects
consist of
sequences of  decreasing subsets of $M$ of the  form $ (s_.) = (M = s_1
\supsetneqq s_2 \ldots \supsetneqq s_p \neq \emptyset),\;\; p>0.$

\noindent In this construction we will need double complexes, more
precisely complexes of the previously defined  exterior complexes,
so we introduce {\it  the category $M_.^+$ whose objects are the
subsets $J\subset M$ including the empty set } so that the DeRham
complex {\it $\Omega (L,  N.)$ is written now as $s(L_J)_{J
\subset M}$} and we consider  objects with indices in the category
$M_.^+\times S(M)$.\\
 Geometrically  $M$ corresponds to a normal section to $Y^*_M$ in $X$
and  $J$  to $\wedge_{i\in J}dz_i$ in the exterior DeRham complex
on the normal section to  $Y^*_M$. The decomposition $M_.^+\simeq
(M-K)_.^+ \times K_.^+$ corresponds to the isomorphism $\C^M\simeq
\C^{(M-K)}\times \C^K$.

\smallskip
\noindent {\it Notations.} For each $s.\in S(M)$ let
$W^{s_{\lambda}} = W (\sum_{i \in  s_{\lambda}}N_i)$ centered at
$0$, for $J\subset M$ and an integer $r$,  we define
$a_{s_{\lambda}}(J,r) = \vert s_{\lambda}\vert  - 2 \vert
s_{\lambda}\cap J\vert + r $,
 and
 for all $(J,s.)\in M_.^+\times S(M)$ the functorial vector
 spaces\\
$ W_r(J,s.)L \colon = \bigcap_{s_{\lambda}\in
s.}W^{s_{\lambda}}_{a_{s_{\lambda}}(J,r)} L,\quad F^r(J,s.)\colon=
F^{r-\vert  J\vert  } L , \quad W^{s_{\lambda}} = W (\sum_{i \in  s_{\lambda}}N_i)$,\\
  then we consider for each $(s.)$ a
DeRham complex $\Omega (L,N.)$

 \smallskip
\noindent  {\bf Definition.} {\em The weight ${\cal W}$ (centered
at zero) and Hodge $F$
 filtrations on  the combinatorial  DeRham complex
 $\Omega^{*} L = s(\Omega (L,N.))_{s.\in S(M)}$ are defined by
 ``summing'' over $J$ and $s.$
 $$(\Omega^{*} L ,{\cal W},F) \leqno {(9)}$$
where $\;{\cal W}_r (\Omega^{*} L)\colon=
s(\bigcap_{s_{\lambda}\in
s.}W^{s_{\lambda}}_{a_{s_{\lambda}}(J,r)} L)_{(J,s.) \in
M_.^+\times S(M)},\; \; a_{s_{\lambda}}(J,r) = \vert
s_{\lambda}\vert  - 2 \vert  s_{\lambda}\cap J\vert + r $\\
and   $F^{r} (\Omega^{*} L) \colon=  s(F^{r-\vert J \vert  } L
)_{(J,s.) \in M_.^+\times S(M)}\,$.}

 \smallskip
\noindent The filtrations can be constructed in two times, first by summing
over $J$ to get the  sub-complexes\\
$   W_r(s.)= s(W_r(J,s.))_{J\subset M} $ (weight) and $
F^r(s.)\colon= s(F^{r}(J,s.))_{J\subset M}$ (Hodge).

 \smallskip
 \noindent {\em  Example in dimension $2$}.

\noindent Let $W^{1,2} = W(N_1 + N_2), W^1 = W(N_1 ) $ and $ W^2 =
W(N_2) $, the weight ${\cal W}_r $ is a double complex:\\ ${\cal
W}_r ( \{1,2 \}\supsetneqq 1 )\oplus {\cal W}_r (\{1,2
\}\supsetneqq 2) \rightarrow {\cal W}_r  (\{1,2 \}) $ \\ where the
first line is the direct
 sum  of:

$ {\cal W}_r ( \{1,2 \}\supsetneqq 1 ) = (W_{r+2}^{1,2}\cap
W^1_{r+1}\stackrel{(N_1,\; N_2)} {\longrightarrow} W_{r}^{1,2}\cap
W^1_{r-1}\oplus W_{r}^{1,2}\cap W^1_{r+1}\stackrel{(- N_2,\;N_1)}
\longrightarrow W_{r-2}^{1,2}\cap W^1_{r-1})$

\noindent and

$ {\cal W}_r ( \{1,2 \}\supsetneqq 2 ) = (W_{r+2}^{1,2}\cap
W^2_{r+1}\stackrel{(N_1,\; N_2)} {\longrightarrow} W_{r}^{1,2}\cap
W^2_{r+1}\oplus W_{r}^{1,2}\cap W^2_{r-1} \stackrel{(- N_2,\;N_1)}
\longrightarrow W_{r-2}^{1,2}\cap W^2_{r-1})$

\noindent The second line for $ \{1,2 \}$ is

$ {\cal W}_r ( \{1,2 \} ) = (W_{r+2}^{1,2} \stackrel{(N_1,\; N_2)}
{\longrightarrow}
 W_{r}^{1,2}\oplus W_{r}^{1,2}\stackrel{(- N_2,\;N_1)} \longrightarrow
W_{r-2}^{1,2})$.

\noindent which reduces to the formula in [26] for  $ r = -1 $.

\be { The Complexes $C_r^{KM}L$ and $ C^K_r L $.}{1.II.6} \ To
study the graded part of the weight, we need to introduce the
following subcategories:

\noindent For each subset $ K \subset M$, let $ S_K(M) = \{ s. \in
S(M): K \in s. \}$ (that is $\exists \lambda : K = s_{\lambda}$).
The isomorphism of categories:
$$S(K)\times S(M-K) {\buildrel \sim \over \rightarrow }
 S_K(M), \; (s., s.') \rightarrow  ( K \cup s.', s.)$$
will be of important use later.
 We consider the vector spaces with indices $(J,s.) \in
M.^+ \times S_KM,$ \\
 $C^{KM}_rL(J,s.)\colon=  \bigcap _{K\neq s_{\lambda}\in s.}
W^{s_{\lambda}}_{a_{\lambda}(J,r-1)}  Gr^{W^K}_{a_{K(J,r)}} L $
and for each $(s.)$ the associated complex obtained by summing
over $J$ (resp. over $(s.)$:\\
$\;C^{KM}_rL(s.) \colon= s(C_rL(J,s.))_{J \in M^+_.}$,
$\;C^{KM}_rL \colon= s(C^{KM}_rL(s.))_{(s.) \in S_KM }$.\\
We write $C^{K}_rL(J,s.)$, $C^{K}_rL(s.)$ and $C^{K}_rL$ when $K =
M$.

\smallskip
\noindent {\bf Definition}. {\em For $K \subseteq  M$ the complex
$C_r^{KM} L $ is defined by summing over $J$ and ($s.$)
$$ \;C_r^{KM}L\colon=  s \, (\bigcap _{K\neq s_{\lambda}\in s.}
W^{s_{\lambda}}_{a_{\lambda}(J,r-1)}  Gr^{W^K}_{a_{K(J,r)}} L) \,
_{(J,s.)\in M.^+ \times S_K(M)}  \leqno{(10)}$$
 In the case  $ K = M $ we write
$C_r^{K} L $  }
 $$C_r^K L = s(C^{K}_rL(J,s.))_{(J,s.) \in K.^+ \times S(K)} = s(
(\cap_{K\supsetneqq s_{\lambda}\in s.}
W^{s_{\lambda}}_{a_{s_{\lambda}}(J,r-1)}Gr^{W^K}_{a_{K(J,r)}}L)_{(J,s.)
\in K_.^+ \times S(K)} \leqno{(11)}$$
 \subsection { Purity of the cohomology of the complex $C^K_r L$ } In
this subsection we aim to prove that the filtration ${\cal W}$
will lead to  the weight of what would be in the proper case
 a mixed Hodge complex in Deligne's terminology, that is the induced
 filtration by $F$ on
 the graded parts $Gr^{\cal W}\Omega^{*} L $  is a Hodge filtration.
 For this we need to
 decompose the complex as a direct sum of intermediate extensions of
 variations of Hodge
 structures (which has a meaning locally) whose cohomologies are  pure Hodge
structures [5] and [24] in the proper case. The decomposition
itself is in  the next section.
 Here we prove  the purity of the
complex $C^K_r L$. Its unique non vanishing cohomology will  the
fiber of the variations of Hodge structures  needed in the
decomposition of $Gr^{\cal W}\Omega^{*} L $. The result here is a
fundamental step in the general proof. The key lemma  proved
earlier provides what seems to be the elementary property
 at the level of a nilpotent orbit that leads to establish
the purity and decomposition results. The proof of the theorem
below will occupy the whole subsection. First we present a set of
elementary complexes. Second we prove the purity result on the
complexes $C^K_r L$ which behave as a direct sum of elementary
complexes.

\smallskip
\noindent   Let $L$ be a polarised nilpotent orbit, then the
complexes $C^K_r L $  satisfy the following properties
 \be {Proposition}{3.4}{\em (Purity). The  cohomology of the complex $C_r^K L$,
concentrated in a unique degree, underlies a polarised  $ HS$}.

\smallskip
\noindent The proof of the proposition will occupy the whole
subsection and is divided in two parts. Precise information can be
found in the proposition below.

 \be {Elementary complexes.}{3.5}
 We suppose $K $ of length $\vert K
\vert  = n $ and we identify $K$ with the set of integers $[1,n]$,
then the  elementary complexes are  defined by the following simplicial
vector spaces. For  $J \subset [1,n], $ let
$$ K((m_1,\cdots,m_n),J) L = Gr^{W^n}_{m_n - 2 \vert \{ n\} \cap J \vert  }
\cdots Gr^{W^i}_{m_i - 2 \vert \{ i\} \cap J \vert  } \cdots
Gr^{W^1}_{m_1 - 2 \vert \{ 1\} \cap J \vert  } L \leqno{(12)}$$
  The endomorphism $ N_i $ induces a morphism
denoted also \\ $ N_i: K((m_1,\cdots,m_n),J) L \rightarrow
K((m_1,\cdots,m_n),J\cup {i})L $ trivial for $i \in J$.\\
{\it Remark}. Instead of $L$, we can consider such formulas for
various natural spaces derived from $L$ such as $L/N_J L $, $ N_J
L $
 or $\cap_{ s_{\lambda}\in s.} W^{s_{\lambda}}_{a_{\lambda}}L$
 for a sequence $s.$ of subsets of $[1,m]$ containing $[1,n]
\subset [1,m]$.

\smallskip
\noindent {\bf Definition}. {\em The elementary complexes are the
simple
 complexes associated to the simplicial vector spaces (12) by summing over
 $J \subseteq [1,n]$}
$$K(m_1,\cdots,m_n) L \, \colon= s( K((m_1,\cdots,m_n),J)
L , N_i )_{J \subseteq [1,n]}\leqno{(13)} $$
  \noindent {\bf Proposition}: {\em
For any $((m_1,\cdots,m_n) \in
\Z^n $ let $J(m.) = \{i \in  [1,n]: m_i > 1\}  $.\\
 The cohomology of an elementary complex
 $K(m_1,\cdots,m_n) L$ is isomorphic to a sub-quotient of the vector space
  $K((m_1,\cdots,m_n), J(m.))L$
 concentrated in degree $\vert J(m.) \vert  $. Moreover it vanishes
if there exists at least one $ m_i = 1$.

\smallskip
\noindent More precisely, the cohomology is isomorphic to\\
  $ \;\;  K((m_1,\cdots,m_n),J(m.) ) [(\cap_{i \notin J(m.)} (ker N_i:
L/(\Sigma_{j \in J(m.)  } N_j L) \rightarrow L/
(\Sigma_{j \in J(m.)} N_j L)] \simeq $\\
$Gr^{W^n}_{m_n - 2 \vert \{ n\} \cap J(m.) \vert  } \cdots
Gr^{W^1}_{m_1 - 2 \vert \{ 1\} \cap J(m.) \vert  } [(\cap_{i
\notin J(m.)} (ker N_i: L/(\Sigma_{j \in J(m.)  } N_j L)
\rightarrow L/(\Sigma_{j \in J(m.)} N_j L))]$.}

\smallskip
 \noindent  The proof by induction on $n$ is based on the
  fact that given  an index $i$, we can view $K(m_1,\cdots,m_n) L$
  as the cone over

 \noindent $N_i:
K(m_1,\cdots,\widehat{m_{ i}},\cdots, m_n) (Gr^{W^i}_{m_i}L)\;
{\buildrel N_i \over \rightarrow }
\;K(m_1,\cdots,\widehat{m_{i}},\cdots, m_n)
(Gr^{W^i}_{m_i-2}L)$.\\
 It is enough to notice that $N_i: Gr^{W^i}_{m_i}L \to
Gr^{W^i}_{m_i-2}L$ is injective if $m_i > 0$, surjective if $m_i <
2$ ( bijective for $m_i =1$). The associated morphism on the
complex $ K(m_1,\cdots,\widehat{m_{i}},\cdots, m_n)$ will have the
same property  since the constituent vector spaces respect exact
sequences by strictness of $MHS$.

\noindent Hence if $m_i > 0 $ (resp. $m_i < 2$), $ N_i $ is
injective on $Gr_{m_i}^{W^i}L $ (resp. surjective onto
$Gr_{m_i-2}^{W^i}L $ ) and

\noindent $ K(m_1,\cdots,m_{i},\cdots, m_n) L
\cong  K(m_1,\cdots,\widehat{ m_{i}},\cdots, m_n)
(Gr^{W^i}_{m_i - 2}(L/N_iL))[-1]$

\noindent (resp. $K(m_1,\cdots,\widehat{m_{i}},\cdots, m_n)
(Gr^{W^i}_{m_i}(ker N_i: L\rightarrow  L)$)

\noindent where $K(m_1,\cdots,\widehat{m_{i}},\cdots, m_n) $ is
applied to the polarised nilpotent orbit $Gr^{W^i}_{m_i-2}( L/N_i
L)$

\noindent (resp.$ Gr^{W^i}_{m_i}(ker N_i: L\rightarrow  L)$) with
the nilpotent endomorphisms $ N_j'$ induced by $N_j$ for $j \neq
i$.\\
{\it Remark}. The cohomology space is symmetric in the operations
kernel and cokernel and is isomorphic to

\noindent $K((m_1,\cdots,m_n),J(m.))  [(\cap_{\{i: m_i < 1 \}} ker
N_i)/(\Sigma_{\{j : m_j > 1 \}} N_j (\cap_{\{i: m_i < 1 \}} ker
N_i))]$.

\noindent that is at each process of taking $Gr_{m_i}^{W^i} $ we
apply
 the functor $ ker $ if $m_i \notin J(m.) $ and $coker $ if $m_i \in
 J(m.)$
 \be {Purity of $C^K_r L$.}{3.5}
 {\it Decomposition into combinatorial elementary complexes.}

\n  By the natural decomposition of the relative filtrations in
the Key lemma,  we have isomorphisms, functorial for the
differentials of $C^K_r L$\\
$$ Gr^{W^K}_{a_{K(J,r)}} ( \cap _{K\supsetneqq s_{\lambda}\in s.}
 W^{s_{\lambda}}_{a_{\lambda}(J,r-1)}) L \simeq
 \oplus_{m. \in X(J,s.,r)} Gr^{W^n}_{m_n}\cdots Gr^{W^i}_{m_i  }
 \cdots Gr^{W^1}_{m_1} L,\;   $$
where for all $ (J,s.)\in K_.^+ \times S(K), $\\
$ X(J,s.,r) = \{ m. \in {\Z}^n :\Sigma_{i \in K} m_i =
a_{K}(J,r)\; \hbox {and} \; \forall s_{\lambda}\in s.,
s_{\lambda}\neq K,  \Sigma_{i \in s_{\lambda } } m_i \leq
a_{\lambda}(J,r-1)  \}$\\
 \noindent In
particular, if we define for $J = \emptyset,
 X(s.,r) =  X(\emptyset,s.,r) $ as
$$ X(s.,r) = \{ m. \in {\Z}^n :\Sigma_{i \in K} m_i =
\vert  K \vert  + r \;\; \hbox {and} \;\; \forall s_{\lambda}\in
s., s_{\lambda}\neq K, \Sigma_{i \in s_{\lambda } } m_i \leq \vert
s_{\lambda}\vert   + r-1  \}$$ \noindent  the complex $C_r^K L
(s.)$ splits as a direct sum of
elementary complexes\\
$ C_r^K L (s.)\simeq \oplus_{ X(s.,r)}  K(m_1,\ldots,m_n) L. $

\smallskip
\n {\it The combinatorial elementary complex}. For each $(m_1,
\ldots, m_n) \in \Z^n$  we define a
 complex with indices in $s. \in S(K)$ as follows: \\
 $K(m_1,\ldots,m_n;r)L(s.)=
 K(m_1,\ldots,m_n) L$ if $(m_1,\ldots,m_n) \in X(s.,r)$ and $0 $
 otherwise.

 \smallskip
\n {\bf Definition}. {\em The combinatorial elementary complex is
defined by summing over $s.$}
 $$\tilde  K(m_1,\ldots,m_n;r)L = s[K(m_1,\ldots,m_n;r)L(s.)]_{s. \in
 S(K)}$$
{\bf Lemma}. {\em Define $ X(r) = \{ m. \in {\Z}^n :\Sigma_{i \in
K} m_i = \vert  K \vert  + r \}$,  then we have the decomposition:
$$C_r^K L \simeq \oplus_{m. \in X(r)} \tilde K(m_1,\ldots,m_n;r).\leqno{(14)}$$}
{\it  Example.} For $K = \{ 1,2 \}$, $s.$ is one of the $3$
elements
  $s.' = \{ 1,2 \} \supset \{ 1 \}$, $s.'' = \{ 1,2 \} \supset \{ 1 \}
$ or $ K = \{ 1,2 \}$. Then $\tilde  K(m_1,\ldots,m_n;r)L$ is defined for fixed
$r$ by the combinatorial complex\\
$K(m_1,m_2;r)L(s.') \simeq  K(m_1, m_2) L$
if $m_1 + m_2= 2+r$ and $ m_1 \leq r$ and $0$ otherwise\\
$K(m_1,m_2;r)L(s.'') \simeq   K(m_1, m_2) L$
if $m_1 + m_2= 2+r$ and $ m_2 \leq r$ and $0$ otherwise\\
$K(m_1,m_2;r)L(K) \simeq  K(m_1, m_2) L$ if $m_1 + m_2= 2+r$  and
$0$ otherwise.

\smallskip
\n Notice $\tilde K(m_1,m_2;r)(W^1_1 L)\cong 0 $ for $r \geq 0$
that is the case $m_1 < 2$. In fact suppose $r>0$ in the example
and $m_1 < 2$, then $m_2 = 2+r-m_1 > r$, hence
$K(m_1,m_2;r)L(s.'') \simeq 0$. Moreover $K(m_1,m_2;r)L(s.')
\simeq K(m_1,m_2;r)L(K)$ so that $\tilde K(m_1,m_2;r)(W^1_1 L)
\cong 0$.\\  This is a  main point that we  prove in a more
general setting in the coming basic lemma.

\smallskip
 \n  { \it The elementary
sub-complexes supporting the cohomology}. \\
 For
$r > 0$ we define \\
$T(r) = \{(m_1,\cdots,m_n) \in {\bf N}^n: \forall i \in K, m_i
\geq 2 $ and $ \Sigma_{i \in K } m_i  = \vert  K \vert  +r \}$ \\
( for $ r = 0, \, T(0) = \emptyset$ ) so to introduce the complex
$$ C(T(r)) L \simeq \oplus_{(m_1,\cdots,m_n) \in T(r)} K(m_1, \ldots,
m_n) L \leqno {}$$
 \noindent Dually, for $r < 0$ we define \\
$T'(r) = \{(m_1,\cdots,m_n) \in {\bf Z}^n: \forall i \in K, m_i
\leq 0 $ and $ \Sigma_{j \in K } m_j  = \vert  K \vert  +r \}$ \\
so to
introduce the complex
$$ C(T'(r)) L \simeq \oplus_{(m_1,\cdots,m_n) \in T'(r)} K(m_1, \ldots,
m_n) L \leqno {}$$
 {\bf Lemma.}{\em  The complex $ C(T(r))L$  embeds diagonally into the direct sum
of $C_r(s.)L$ for all sequences $(s.)$ of maximal length $\vert K
\vert $ so to define a morphism of complexes $ C(T(r))L \to C^K_r
L
$.\\
 Dually, the complex $C(T'(r))$ embeds in $C_r(s.)L$ for $ s.
= K$ consisting of one subset $K$, so to define a morphism of
complexes $ C(T'(r))L [1- \vert K \vert ]\to C^K_r L $.}

\n Proof. We check the conditions defined by  all $ s_{\lambda}
\in (s.),\, s_{\lambda}\neq K $, namely $ \Sigma_{i \in s_{\lambda
} } m_i \leq \vert s_{\lambda}\vert  + r-1 $ for all $(m.) \in
T(r)$ by induction on $ \vert  s_{\lambda}\vert  $. We start with
the condition  $\Sigma_{j \in K } m_j = \vert  K \vert +r $
defined by $K$, which is satisfied by definition of $T(r)$. Let $
k \in K$, then $ \Sigma_{j \in K - \{k \}} m_j \leq \vert  K - \{k
\} \vert + r-1 $ by subtracting $m_k
> 1$, which proves the assertion for $s_{\lambda} = K - \{k \}$, hence for
all $s_{\lambda}$ such $\vert  s_{\lambda}\vert  = \vert  K \vert
-1$. Let $A \subset K, A \neq K$ and suppose $\Sigma_{j \in A }
m_j \leq \vert  A \vert + r-1 $ true, by induction on $\vert  A
\vert$, then for $k \in A, B = A - \{ k \}$ we deduce $\Sigma_{j
\in B } m_j \leq \vert  B \vert + r - 2 $ by subtracting $m_k
> 1$.\\
Dually, for $(s.) = K$  there is no additional conditions, so  the
statement is clear.
  \be {Proposition}{3.6}
\n i){\em For  $r > 0$, the canonical embedding of $ C(T(r)) L$
 into  $C_r^K L$ induces an isomorphism on the
cohomology.} \\
In particular the cohomology of $C_r^K L$, concentrated in degree
$\vert  K \vert $, is isomorphic to

 $ H^{\vert  K \vert}(C_r^K L) \simeq Gr^{W^K}_{r- \vert  K \vert }
 [L/(\Sigma_{i \in K} N_i
L)] \simeq \oplus_{(m.) \in T(r)} Gr^{W^n}_{m_n - 2 } \cdots
Gr^{W^1}_{m_1 - 2 } [L/(\Sigma_{i \in K} N_i L)] $

\noindent it is a polarised $HS$ of  weight $r + m - \vert  K
\vert $ with the weight filtration induced by $W^K$ shifted by
$m$ and Hodge filtration induced by $F^K$.\\
ii) If $r = 0$, the complex $C_r^K L$ is acyclic.

\n iii){\em Dually, for $ r <0 $, the canonical embedding of $
C(T'(r))L [1 - \vert  K \vert ]$
 into  $C_r^K L$ induces an isomorphism on the
cohomology. \\
  In particular the cohomology of $C_r^K L$, concentrated in
  degree $\vert  K \vert  - 1$,
  is isomorphic to

\noindent $ H^{\vert  K \vert - 1}(C_r^K L) \simeq
Gr^{W^K}_{r+\vert K \vert } (\cap_{i \in K} (ker N_i: L \to L))
\simeq \oplus_{(m.) \in T'(r)} Gr^{W^n}_{m_n } \cdots
Gr^{W^1}_{m_1 } [\cap_{i \in K} (ker N_i: L \to L )]$

\noindent it is a polarised $HS$ of  weight $r + m + \vert  K
\vert $ with the weight filtration induced by $W^K$ shifted by $m$
and Hodge filtration induced by $F^K$}.

\smallskip
\noindent {\em Remark}:  If $ r \in [1, \vert  K \vert  - 1]$,
$T(r)$ is empty and $C^K_r L$ is acyclic.  If $ r \in [ - \vert  K
\vert  + 1, 0]$,
 $T'(r)$ is empty and $C^K_r L$
is acyclic. In all cases $C^K_r L$ appears in $Gr^{{\cal W}}
{\Omega}^* L$.

\smallskip
\noindent The principal ingredient in the proof is based on:\\
{\bf Lemma (basic)}. {\em i) For $ r \geq 0$, the complex $\tilde
K(m_1,\ldots,m_n;r)L$
  is acyclic whenever
  at least one $m_i < 2$. Equivalently, for each  $i \in K$,
  the complex $C_r^K (W_{1}^iL)$
  is acyclic.\\
  ii) Dually, for  $r \leq 0$, the complex $\tilde  K(m_1,\ldots,m_n;r)L$
  is acyclic whenever
  at least one $m_i \geq  2$. Equivalently, for each  $i \in K$, the complex
  $C_r^K (L/W_{0}^i L)$ is acyclic.}

\smallskip
\noindent  The equivalences follow from the decompositions \\
  $ C_r^K (W_{1}^iL) \simeq \oplus_{m. \in X(r), m_i < 2} \tilde
  K(m_1,\ldots,m_n)L$\\
   $ C_r^K (L/W_{1}^iL) \simeq \oplus_{m. \in X(r), m_i \geq 2}
  \tilde K(m_1,\ldots,m_n,r)L$.\\
  We note that $ \tilde K(m_1,\ldots,m_n,r)L \cong 0$ if at least
  one $m_i= 1$.

\smallskip
\noindent  Proof of the lemma.  The result is based on the
following elementary remark. Consider in $s.$ a  sequence $s_{a+2}
=  s_{a+1} \cup \{i\} \supset s_{a+1}$ with $i \notin s_{a+1}$,
then the condition on $m.\in X(s.,r)$ associated to $s.$ defined
by $s_{a+1} \cup \{i\}$   is $m_i + \Sigma_{j \in s_{a+1}} m_j
\leq \vert s_{a+1} \cup \{i\} \vert + r-1 $ to compare with the
condition $\Sigma_{j \in s'_{a+1}} m_j \leq \vert s'_{a+1} \vert +
r-1 $ defined by  $ s_{a+1}$. Precisely when $m_i < 2$ (that is in
$W^i_1 L $) the condition for $ s_{a+1}\cup \{i\}$ follows from
the condition for $s_{a+1}$, hence the conditions defined by the
subsequence $s_{a+1} \cup \{i\} \supset s'{a+1}$ in $s.$ is the
same as the condition defined by $ s_{a+1}$ in $ d_{s_{a+2}}(s.)$
where $ d_{s_{a+2}}$ is the differential consisting in the removal
of $s_{a+2}$, so that in the sum over all $s.$ this couple is
quasi-isomorphic to zero. The following rigorous proof consists on
filtering the complex by carefully choosing subsets in $S(K)$. All
constructions below are compatible with the decompositions
 and apply
to each  complex  $\tilde K(m_1,\ldots,m_n,r)L $.\\
 i) We construct a filtration of $C^K_r (W_1^i L)$ with
acyclic sub-complexes. We define the $i$-length $\vert s.
\vert_{i} $ of a sequence $s.$   as the number of subset
$s_{\lambda}$ not containing $i$. Let $S_a $ be the full
subcategory whose objects satisfy $\vert  s. \vert_{i} \leq a $,
hence an object $s.$ in $S_a $ is written as $s.''\cup \{i\}
\supset s.'$ with $\vert s.' \vert \leq a $ with no subset in $
s.'$  containing $i$. Deleting a subset of $s.$ in $S_a $ gives
another object in
 $S_a $, hence $C(S_a)L \colon=  s[C_r^K L (s.)]_{s. \in S_a }$
is a sub-complex of $C_r^K L$ where $C_r^K L (s.)$ is obtained by
summing over
over $J \subset K$ for a fixed $(s.)$.\\
1) To start we write $C(S_0)L $ as a cone  over a morphism
inducing a quasi-isomorphism on $C(S_0)( W^i_1L) $. We divide the
objects of $S_0$ in two families : $S'_0 $ whose objects are
defined by the sequences $(s.)$ starting with the subset $\{i\}$
and $S''_0 $ whose objects are defined by the sequences $(s.)$
whose elements $s_{\lambda}$ contain $i$ but are different from $
\{i\}$. Then we consider the two complexes $C(S'_0)L \colon=
s[C_r^K L(s.)]_{s. \in S'_0 }$ and $C(S''_0)L \colon=
s[C_r^K L(s.)]_{s. \in S''_0 }$.\\
 Let $(s.) = K \supset \ldots
s_{\lambda}\supset \ldots \supset \{i\}$ in $ S'_0$; deleting
$s_{\lambda}\neq \{i\} $ is a morphism in $S'_0$ but deleting
$\{i\}$ gives an element in $S''_0$ . It is easy to check that
deleting $\{i\}$ defines a morphism of complexes $ I_{\{i\}}:
C(S'_0 )L \rightarrow C(S''_0)L $ defined by embedding $C_r^K
L(s.)$ into $C_r^K L(d_{\{i\}} (s.))$ where $d_{\{i\}} (s.)= (s.)
- \{i\}$. The cone over this $ I_{ \{i\} }$ is isomorphic to
$C(S_0)L[1] $.

\noindent {\it We show now that if we reduce the construction to
$W_1^iL$ instead of $L$, the morphism $ I_{\{i\}  }: C(S'_0
)(W_1^iL) \rightarrow C(S''_0)(W_1^iL) $ is an isomorphism.} Let
$s. = (K \supset \ldots s_{\lambda}\supset \ldots \supset \{i\})$
in $ S'_0$. It is enough to notice that the condition on $m. \in
X(s.,r)$ associated to $\{i\}$ defined by $m_i \leq \vert \{i\}
\vert + r-1 = r $ is irrelevant since $m_i \leq 1$ and $r > 0$ (
since for $m_i = 1$ the elementary complexes are acyclic, we could
also suppose $m_i < 1 $ and $ r = 0 $) . It follows that
$C(S_0)(W_1^iL)$ is acyclic.

\smallskip
\noindent 2) {\it We extend the proof from $S_0$ to $S_a$}.\\
Suppose by induction that $C(S_a)(W_1^iL) $ is acyclic for $a \geq
0$, we prove $C(S_{a+1})(W_1^iL)$ is also acyclic. It is enough to
prove that the
quotient $G_a L \colon= C(S_{a+1})L/ C(S_{a})L$ is acyclic.\\
Let $s. = s.''\cup \{i\} \supset s.'$ with  $\vert  s.'  \vert
\leq a+1 $ ( $ s.'$ not containing $i$).  Deleting $s.''\cup \{i\}
$ is a morphism in $S_{a+1}$ but deleting $s.'$ gives an element
in $S_{a}$, hence defines a differential zero in $G_a L$.\\
 We divide the objects of $S_{a+1}$ in two families containing $S_a$
 according to the
subsequence $s''_{a+2} \cup \{i\} \supset s'_{a+1}$, the family
$S'_{a+1} $ (and respectively  $S''_{a+1} $) whose objects are
defined by the sequences
 satisfying $s''_{a+2} = s'_{a+1}$
(resp. $s''_{a+2} \varsupsetneq s'_{a+1}$).  Deleting
$s''_{\lambda} \cup \{i\} $ for
 $s''_{\lambda} \neq s''_{a+2} $ is a morphism in $S'_{a+1}$ but for
 $s''_{\lambda} = s''_{a+2} $ it defines a functor $d_{s_{a+2}}:
S'_{a+1} \rightarrow S''_{a+1}$. If we consider the complexes
$C(S'_{a+1})L/C(S_a)L$ and $C(S''_{a+1})L/C(S_a)L$, we deduce a
morphism of complexes $ I_{s_{a+2}}:  C(S'_{a+1})L/C(S_a)L
\rightarrow C(S''_{a+1})L/C(S_a)L$ which consists in an embedding
of $C_r^K L(s.)$ into $C_r^K L(d_{s_{a+2}} (s.))$ where
$d_{s_{a+2}} (s.) = (s.)-s_{a+2}$. It is easy to check that the
cone over this $ I_{s_{a+2}}$ is isomorphic to $ G_a L[1] $.

\smallskip
\noindent {\it  We show now that if we reduce the construction to
$W_1^iL$ the morphism $ I_{s_{a+2}}$ is an isomorphism.} The
condition on $m.\in  X(s.,r)$ associated to $s.$ defined by
$s''_{a+2} \cup \{i\}$ when $ s''_{a+2} = s'_{a+1}$  is   $m_i +
\Sigma_{j \in s'_{a+1}} m_j \leq \vert s'_{a+1} \cup \{i\}  \vert
+ r-1 $ to compare with the condition  $\Sigma_{j \in s'_{a+1}}
m_j \leq \vert s'_{a+1} \vert + r-1 $ defined by  $ s'_{a+1}$).
Precisely when $m_i < 2$ (that is in $W^i_1 L $) the condition for
$ s'_{a+1}\cup \{i\}$ follows from the condition for $s'_{a+1}$,
hence the conditions defined by the subsequence $s'_{a+1} \cup
\{i\} \supset s'{a+1}$ in $s.$ is the same as the condition
defined by $ s'_{a+1}$ in $ d_{s_{a+2}}(s.)$. This proves that $
I_{s_{a+2}}$ is an isomorphism for $W_1^iL$, hence $G_a (W_1^iL)$
is acyclic and i) follows by induction.

\smallskip
\noindent  ii) {\it Dual proof} (to be skipped). We construct a
dual filtration of $C^K_r (L/W_0^i L)$ with acyclic sub-complexes.
To simplify notations we denote a sequence by $s^.$ and define its
$i$-colength  $\vert s^. \vert^i $ as the number of subsets
$s_{\lambda}$ containing $i$. Let $S^a $ be the full subcategory
whose objects satisfy $\vert s^.  \vert ^i \leq a $, hence an
object  $s^.$ in $S^a $ is written as $s.''\cup \{i\} \supset
s'^.$ with  $\vert s''^. \vert \leq a $ and $ s'^.$ not containing
$i$. Deleting a subset of $s^.$ in $S^a $ gives another object in
 $S^a $, hence $C(S^a)L \colon=  s[C_r^K(s^.)L]_{s^. \in S^a }$
is a sub-complex of $C_r^K L$.\\
To start with  $C(S^0)L = 0$ since $S^0$ is empty.\\
Suppose by induction that $C(S^a)(L/W_0^iL) $ is acyclic for $a
\geq 0$, we prove $C(S^{a+1})(L/W_0^iL)$ is also acyclic. It is
enough to prove that the
quotient $G^a L  \colon= C(S^{a+1})L/ C(S^{a})L$ is acyclic.\\
Let $s^. = s''^.\cup \{i\} \supset s'^.$ with  $\vert  s''^. \vert
\leq a+1 $ ( $ s''^.$ not containing $i$). Deleting $s'^.$  is a
morphism in $S^{a+1}$ but deleting $s''^.\cup \{i\} $ gives an
element in $S^{a}$, hence defines a differential zero in $G^a
L$.\\
 We divide the objects of $S^{a+1}$ in two families containing $S^a$
 according to the subsequence $s''^{a+1} \cup \{i\} \supset
s'^{a+2}$: $S_1^{a+1} $ (resp. $S_2^{a+1} $) whose objects are
defined by the sequences  $s''^{a+1} = s'^{a+2}$ (resp. $s''^{a+1}
\varsupsetneq s'{a+2}$).  Deleting $s'^{\lambda}  $ for
 $s'^{\lambda} \neq s'^{a+2} $ is a morphism in $S_1^{a+1}$ but for
 $s'^{\lambda} = s'^{a+2} $ it defines a functor $d_{s^{a+2}}:
S_1^{a+1} \rightarrow S_2^{a+1}$. If we consider the complexes
$C(S_1^{a+1})L/C(S^a)L$ and $C(S_2^{a+1})L/C(S^a)L$, we deduce a
morphism of complexes $ I_{s^{a+2}}:  C(S_1^{a+1})L/C(S^a)L
\rightarrow C(S_2^{a+1})L/C(S^a)L$ by embedding $C_r^K L(s^.)$
into $C_r^K L(d_{s^{a+2}} (s^.))$. It is easy to check that the
cone over this $ I_{s^{a+2}}$ is isomorphic to $ G^a L [1] $.

\smallskip
\noindent We show now that if we reduce the construction to
$L/W_0^iL$ the morphism $ I_{s^{a+2}}$ is an isomorphism. The
condition on $m.\in  X(s^.,r)$  defined by $s''^{a+1} \cup \{i\}$
when $ s''^{a+1} = s'^{a+2}$ (resp. by $ s'^{a+2}$) is $m_i +
\Sigma_{j \in s'^{a+2}} m_j \leq \vert s'^{a+2} \cup \{i\} \vert +
r-1 $ (resp. $\Sigma_{j \in s'^{a+2}} m_j \leq \vert s'^{a+2}
\vert + r-1 $), but precisely when $m_i > 0$ (that is in $L/W^i_1
L $) the condition for $ s'^{a+2}$ follows from the condition for
$s'^{a+2}\cup \{i\}$, hence the conditions defined by the
subsequence $s'^{a+2} \cup \{i\} \supset s'^{a+2}$ in $s^.$ is the
same as the condition defined by $s'^{a+2} \cup \{i\}$ in $
d_{s^{a+2}}(s^.)$. This proves that $ I_{s^{a+2}}$ is an
isomorphism for $L/W_0^iL$, hence $G^a(L/ W_0^iL)$ is acyclic.\\
This shows  by induction that $C(S^{\vert K \vert - 1
})(L/W_0^i)L$ is acyclic. At the last step, we show that $ C^K_r
(L/W_0^iL)$ is quasi-isomorphic to $ C(S^{\vert K \vert
})(L/W_0^i)L = C(S^{\vert K \vert - 1 })(L/W_0^i)L$. In colength
$\vert K \vert$ there is only one sequence not in $ S^{\vert K
\vert - 1 }$ that is the full sequence starting with $\{i \}$
which imposes the condition $m_i \leq r$ which is impossible since
$r \leq 0 $ and $m_i \geq 1$. This ends the proof of the lemma.
 The proposition follows immediately from

\smallskip
\noindent {\bf Corollary.} {\em i) For $r > 0 $, the complex
$C(T(r))$ is contained in each $C^K_r L(s.)$ that is $T(r)\subset
X(s.,r)$ and the complex $\tilde C(T(r))= s(C(T(r))_{s. \in S(K)}$
is
contained in and  quasi-isomorphic to $ C^K_r L$ . \\
ii) Dually, the complex $C(T'(r))$ is contained only in $C^K_r
L(s.)$ for $s. = K$ and is quasi-isomorphic to $ C^K_r L[\vert  K
\vert   - 1]$.}

\smallskip
\noindent  We did check that the complex $C(T(r))$ is contained in
each $C^K_r(s.)$ for $s.$ of maximal length, hence for all $(s.)$.
The lemma shows that $ \tilde K(m_1,\ldots,m_n,r)L \cong 0$
whenever at least one $m_i < 2$, hence i) follows.\\
Dually, the condition for $K$, $\Sigma_{j \in K } m_j  = \vert  K
\vert  +r \Rightarrow \forall A = K - k \subset K, \Sigma_{j \in A
} m_j > \vert  A\vert   + r-1  $ by subtracting   $m_k < 1$. If
this is true for all $A : \vert   A \vert   = a $ then $ \forall B
= (A - k) \subset A, \Sigma_{j \in B } m_j > \vert  B \vert   + r-
1 $ as well. The shift in degree corresponds to the shift for $s.
= K $ in the total complex $ C^K_r L$.
\subsection {Local decomposition.}
Since the  purity result is established, we can easily prove now
the decomposition theorem after a careful study of the category of
indices $S(I)$.
 \be{Theorem}{3.7} {\em (decomposition). For a
nilpotent orbit $L, N_i, i\in M $, $\vert M \vert = n$ and for all
subsets $K \subset M $ there exist canonical morphisms of
$C_r^{KM}L$ in $Gr_r^{{\cal W}} \Omega^{*} L$ inducing a
quasi-isomorphism (decomposition as a direct sum)
$$ Gr_r^{{\cal W}} (\Omega^{*} L) \cong {\oplus_{K \subset M}} C_r^{KM}L.$$
Moreover $ Gr_0^{\cal W} \Omega^{*} L \cong 0  $ is acyclic.}

\smallskip
For $n = 1$, $K$ and $M$ reduces to one element $1$ and the
theorem reduces to

\smallskip
\n \centerline {$Gr_r^{{\cal W}} (\Omega^{*} L) \simeq C^1_r L $
is the complex $
  Gr^{W^1}_{r+1}L {\buildrel N_1 \over
\rightarrow } Gr^{W^1}_{r-1}L $}

\smallskip
\n By the elementary properties of the weight filtration of $N_1$,
it is quasi-isomorphic to $Gr^{W^1}_{r-1} (L/N_1 L) [-1]$ if $r >
0$, $Gr^{W^1}_{r+1}(ker N_1:L \rightarrow L)$ if $r < 0 $ and
$C^1_0 L \simeq  0.$

\smallskip
\noindent The proof is by induction on $n$; we use only the
property $Gr_0^{\cal W} (\Omega^{*} L) \simeq 0 $ in $n-1$
variables to get the decomposition for $ n$, then we use the fact
that $C_0^K L$ for all $K$ is acyclic  to get again $Gr_0^{\cal W}
(\Omega^{*} L) \cong 0 $ for $n$ variables so to complete the
induction step.

\noindent The proof of the decomposition is carried in the three
lemmas below.

\smallskip
\noindent For each $i  \in ${\bf N} we define a map into the
subsets of $M$\\
$ \varphi_i: S(M) \rightarrow {\cal P}( M)$ :\  $ \varphi_i (s.) =
Sup \{s_{\lambda }: \vert  s_{\lambda } \vert \leq i \}$. For each
 $(J,s.) \in M.^+ \times S(M)$, we consider the subspaces
  of $L$ with indices $i$ and $t$\\
 $W_t(i,J,s.)\colon=  (\bigcap _{\varphi_i (s.)
\varsubsetneq s_{\lambda}\in s.}
W^{s_{\lambda}}_{a_{s_{\lambda}}(J,r+ t)}) \cap( \bigcap
_{s_{\lambda} \subset \varphi_i (s.) , s_{\lambda}\in s.}
W^{s_{\lambda}}_{a_{s_{\lambda}}(J,r-1)} L)$

\noindent We define \\
 $ G_i(J,s.)(L) = Gr^{W(i,J,s.)}_0 L = W_0(i,J,s.)/W_{-1}(i,J,s.)$
 and  the complexes
$$G_i(s.)L = s(G_i(J,s.)L)_{J\subset M},\quad
 G_i^{SM} L = s(G_i (s.)L)_{s. \in
S(M)}$$ In particular,  $\varphi_0 (s.) = \emptyset $,
  and $\varphi_{ \vert   M \vert } = M $ so that
 $$G_0(s.)L = Gr^{W(s.)}_r L, \ G_0^{SM}L = Gr^{{\cal W}}_r \Omega^{*} L,
 \quad G_{\vert  M \vert }^{SM}L = C^M_r L$$
 The proof of the  theorem   by induction on
 $i$, starting with $i =0$, is based on

\smallskip
\noindent {\bf Lemma 1. }
$$G_i^{SM}L \cong G_{i+1}^{SM} L \oplus
(\oplus_{ K \subset M, \vert  K \vert = i+1 } C_r^{KM}L)$$ To
relate $G_i^{SM} L$ and $G_{i+1}^{SM} L$, we define $ S^{i}(M) =
\{s. \in S(M): \vert \varphi_i (s.) \vert  \ = i \}$ in $S(M)$\\
and consider the subcategory $S(M) - S^{i+1}(M)$.
 The restrictions of the simplicial  vector spaces $G_i$ and
$G_{i+1}$ to $S(M) - S^{i+1}(M)$, define two
sub-complexes: \\
$G'_i L = s(G_i (s.)L)_{s. \in S(M)- S^{i+1}(M)}$ embedded in $
G_i^{SM} L $ and \\
$ G''_{i+1}L = s(G_{i+1} (s.)L)_{s. \in S(M)- S^{i+1}(M)}$
embedded in $G_{i+1}^{SM} L $\\ since deleting an object  $s. \in
S(M) - S^{i+1}(M)$ gives an object in the same subcategory. We
have $G''_{i+1}L = G'_iL $ since $\varphi_i = \varphi_{i+1}$ on
$S(M) - S^{i+1}(M)$; hence we are reduced to relate the the
quotient complexes: $G_i^{SM} L/G'_i L $ and $G_{i+1}^{SM}
L/G''_{i+1}L$ which are obtained by summing over $S^{i+1}(M)$.\\
 We remark that
  $ S^{i+1}(M)\simeq \cup_{\vert  K \vert  = i+1}S_K(M)$ is a disjoint union of
 $S_K(M)$ where
$\varphi_{i+1} (s.)= K$  with $\vert  K \vert =
i+1$.\\
{\it Definition of $ B_r^{KM} L$.}\\
 For $K$  fixed in $M$, we introduce the subspaces with index $t$\\
 $W_t(K,J,s.)L \colon=  (\bigcap _{K
\varsubsetneq s_{\lambda}\in s.}
W^{s_{\lambda}}_{a_{s_{\lambda}}(J,r+ t)}) \cap (
\bigcap_{s_{\lambda} \subset K , s_{\lambda}\in s.}
W^{s_{\lambda}}_{a_{s_{\lambda}}(J,r-1)} L)$\\
 $Gr^{W(K,J,s.)}_0 L = W_0(K,J,s.)L/W_{-1}(K,J,s.)L$ \\
$B_r^{KM} L = s[Gr^{W(K,J,s.)}_0 L]_{(J,s.)\in M.^+\times
S_K(M)}$\\
It follows by construction\\
 $G_{i+1}^{S^{i+1}(M)}L \colon= s(G_{i+1} (s.)L)_{s. \in S^{i+1}(M)}\simeq
\oplus_{\vert K \vert  = i+1} B_r^{KM} L$,\\ and a triangle
$$ G''_{i+1}L \rightarrow G_{i+1}^{SM} L \rightarrow \oplus_{\vert K \vert  =
i+1} B_r^{KM} L$$
 {\it Definition of $ A_r^{KM} L$}.\\
$ A_r^{KM}(J,s.)L = [(\bigcap _{K \subset s_{\lambda}\in s.}
W^{s_{\lambda}}_{a_{s_{\lambda}}(J,r)}L) \cap
 ( \bigcap_{s_{\lambda} \varsubsetneq K , s_{\lambda}\in s.}
W^{s_{\lambda}}_{a_{s_{\lambda}}(J,r-1)} L)]\;/\;[( \bigcap_{
s_{\lambda}\in s.} W^{s_{\lambda}}_{a_{s_{\lambda}}(J,r-1)} L)]$\\
$ A_r^{KM}L = s (A_r^{KM}(J,s.))_{(J,s.)\in M.^+\times S_K(M)}  $\\
  We
have by construction: $ s(G_i(s.)L)_{s.\in S_KM} = A_r^{KM} L$.
Suppose $\vert K \vert  = i+1$, then \\
 $G_i^{S^{i+1}(M)}L  =
s(G_i (s.)L)_{s. \in S^{i+1}(M)}\simeq
\oplus_{\vert  K \vert  = i+1} A_r^{KM} L $.\\
Considering the triangle
$$  G'_i L \rightarrow  G_i^{SM} L  \rightarrow
 \oplus_{\vert  K \vert  = i+1} A_r^{KM} L$$
and the morphism of triangles: $G_{i+1}^{SM}L \rightarrow
G_{i}^{SM}L $,
  $ B_r^{KM}L \rightarrow A_r^{KM} L $, $ G''_{i+1} L \simeq G'_i
  L$, the relation in the lemma follows from a relation between
$C_r^{KM}L $ and $ A_r^{KM}L $

\smallskip
\noindent { \bf Lemma 2}:
 {\em We have a quasi-isomorphism in the derived category:}
  $$ C_r^{KM}L \oplus B_r ^{KM} L \cong A_r^{KM}L $$
The proof is based on the following elementary remark:

\smallskip
\noindent {\bf Remark.}{\em  Let $W^i$ for $ i =1, 2$  be two
increasing filtrations on an object $V$ of an abelian category and
$a_i$ two integers, then we have an exact sequence:}

\medskip $\quad 0\rightarrow W_{a_2-1}^2 Gr^{W^1}_{a_1}\oplus
W^1_{a_1-1}Gr^{W^2}_{a_2} \rightarrow W^1_{a_1}\cap W_{a_2}^2
/W^1_{a_1-1}\cap W_{a_2-1}^2 \rightarrow
Gr^{W^1}_{a_1}Gr^{W^2}_{a_2}\rightarrow 0 $

\medskip
\noindent We apply this remark to the space \\
 $V = V(J,s.) =
\bigcap _{s_{\lambda} \varsubsetneq  K, s_{\lambda}\in (s.)}
W^{s_{\lambda}}_{a_{s_{\lambda}}(J,r-1)} L,$ \\
filtered by:\\
$W^1_t(K,J,s.)L = W^K_{a_K(J,r+ t)}L \cap V(J,s.)) ,\;\;
 W^2_t(K,J,s.)L = (\bigcap _{K\varsubsetneq s_{\lambda}\in (s.)}
W^{s_{\lambda}}_{a_{s_{\lambda}}(J,r+ t)})L \cap V(J,s.) $ \\
so that \\
$W^1_t(K,J,s.)L \cap W^2_t(K,J,s.) L = (\bigcap _{K\subset
s_{\lambda}\in s.} W^{s_{\lambda}}_{a_{s_{\lambda}}(J,r+ t)})L)
\cap (\bigcap _{s_{\lambda} \varsubsetneq  K, s_{\lambda}\in (s.)}
W^{s_{\lambda}}_{a_{s_{\lambda}}(J,r-1)} L)$.\\
 Let $a_1 = a_2 = t = 0$. , then we deduce  from the above sequence an exact
sequence of vector spaces

$0 \rightarrow     W^2_{-1}(J,s.) Gr^{W^K}_{a_K(J,r)}( V(J,s.))
\oplus W^K_{a_K(J,r-1)}Gr^{W^2}_0 ( V(J,s.))$

$ \rightarrow W^K_{a_K(J,r)}\cap W_0^2 ( V(J,s.)
)/W^K_{a_K(J,r-1)}\cap W_{-1}^2 ( V(J,s.)) \rightarrow
Gr^{W^K}_{a_K(J,r)} Gr^{W^2}_0 ( V(J,s.)) \rightarrow 0 $\\
By summing over $(J,s.)$ in
 $M.^+\times S_K(M)  $ we get an exact sequence of complexes
$$0\rightarrow  C_r^{KM}L \oplus B_r^{KM} L \rightarrow
A_r^{KM}L \rightarrow D_r^{KM}L \rightarrow  0$$ where by
definition
\\
$ D_r^{KM}L = s[Gr^{W^2(J,s.)}_0 Gr^{W^K}_{a_K(J,r)} ( \bigcap
_{s_{\lambda} \varsubsetneq  K, s_{\lambda}\in s.}
W^{s_{\lambda}}_{a_{s_{\lambda}}(J,r-1)} L)]_{(J,s.)\in
 M.^+\times S_K(M)}  $

\noindent Lemma $2$ follows if we prove

\smallskip
\noindent { \bf Lemma 3}. {\em $D_r^{KM}L \cong 0.$}

\noindent  Proof. The idea of the proof is to consider $C^K_r L$
as the fiber of a local system on $Y_K^*$ and form the filtered
complex $(\Omega^{*} (C^K_r L ),{\cal W}) $ for the polarised
nilpotent orbit ($ C^K_r L, N_i, i \in M-K$  to which we can apply
the theorem on the lower dimensional space $Y_K^*$. It happens
that $D_r^{KM}$ is quasi-isomorphic to $Gr^{\cal W}_{0}
(\Omega^{*} (C^K_r L ))$ hence acyclic by induction on dimension.
We can either use that $ C^K_r L $ is reduced to its unique non
zero cohomology or as well
  prove the acyclicity for each term in $C^K_r L $, as we do
  now.

\smallskip
\noindent We consider the  filtration induced by  $W^2$ on
$Gr^{W^K}_{a_K(J,r)}( \bigcap _{s_{\lambda} \varsubsetneq  K,
s_{\lambda}\in s.}W^{s_{\lambda}}_{a_{s_{\lambda}}(J,r-1)} L) $\\
 $W^2_t(K,J,s.)\colon=  \bigcap _{K\varsubsetneq
s_{\lambda}\in s.} W^{s_{\lambda}}_{a_{s_{\lambda}}(J,r+
t)}Gr^{W^K}_{a_K(J,r)} ( \bigcap _{s_{\lambda} \varsubsetneq  K,
s_{\lambda}\in s.} W^{s_{\lambda}}_{a_{s_{\lambda}}(J,r-1)} L)$.\\
 {\it We want to  rewrite the complex  $ D^{KM}_r L$ as a sum in two
times over $(J,s.)\in K.^+\times S(K)$ and $(J',s.')\in
(M-K).^+\times S(M-K))$  corresponding to
$((J \cap J'), (K \cup s.', s.))\in M.^+\times S(M)$}.\\
Using the expression of sequences  in $S_K(M)$  as $ s. = (K \cup
s.',s.)$
 and the result on relative weight filtrations  with respect to
$Gr^{W^K}_{a_K(J,r)}$ we  rewrite $W^2_t$ as \\
$ W^2_t(K,(J,J'), (K\cup s.',s.)) = \bigcap _{s'_{\alpha }\in s.'}
W^{K\cup s'_{\alpha }}_{a_{s'_{\alpha
}}(J',a_K(J,r)+t)}Gr^{W^K}_{a_K(J,r)} ( \bigcap _{s_{\lambda}
\varsubsetneq  K, s_{\lambda}\in s.}
W^{s_{\lambda}}_{a_{s_{\lambda}}((J,J'),r-1)} L)$

\noindent $= \bigcap _{s'_{\alpha }\in s.'}  W^{s'_{\alpha
}}_{a_{s'_{\alpha }}(J',t)}Gr^{W^K}_{a_K(J,r)} ( \bigcap
_{s_{\lambda} \varsubsetneq  K, s_{\lambda}\in s.}
W^{s_{\lambda}}_{a_{\lambda}(J,r-1)} L)$\\
 For a fixed $ (J,s.)\in K.^+ \times S(K)$ we introduce the filtration
  $W_t'^2 = \bigcap _{s'_{\alpha }\in
s.'} W^{s'_{\alpha }}_{a_{s'_{\alpha }}(J',t)}$ on the space
$L(r,J,s.) \colon= Gr^{W^K}_{a_K(J,r)} ( \bigcap _{s_{\lambda}
\varsubsetneq  K, s_{\lambda}\in s.}
W^{s_{\lambda}}_{a_{s_{\lambda}}(J,r-1)} L)$\\
and the complex\\
\noindent $ D(M-K)(L(r,J,s.) \colon= s[Gr^{W'^2}_0(L(r,J,s.))
]_{(J',s.')\in (M-K).^+\times S(M-K)}$\\
\noindent We have by construction
$$ D^{KM}_r L = s[D(M-K)(L(r,J,s.)]_{(J,s.)\in K.^+\times S(K)}$$
We prove by induction on $n$ that each complex
$$ D(M-K)(L(r,J,s.) \colon= s[Gr^{W'^2}_0( Gr^{W^K}_{a_K(J,r)} ( \bigcap _{s_{\lambda}
\varsubsetneq  K, s_{\lambda}\in s.}
W^{s_{\lambda}}_{a_{s_{\lambda}}(J,r-1)} L) ]_{(J',s.')\in
(M-K).^+\times S(M-K)}\cong 0 $$ is acyclic.  Fixing $(J,s.)$, we
decompose $L(r,J,s.)$ into a direct sum of \\
$ L(i.) =
Gr^{W^{i_s}}_{m_{i_s} - 2 \vert {i_s} \cap J \vert } \cdots
Gr^{W^{i_1}}_{m_{i_1} - 2 \vert {i_1} \cap J \vert } L$ for
$\{{i_1},\ldots, {i_s}\} =
K$ and $\sum_{i_p \in K}  m_{i_p} = a_K(J,r)$.\\
We reduce the proof to $ D(M-K)(L(i.)) \cong 0$, then we introduce
the weight filtration ${\cal W}$ on the combinatorial DeRham
complex $\Omega^{*} (L(i.))  $ for the nilpotent orbit $L(i.) $ of
dimension strictly less then $n$ and weight $a_K(J,r) $ and we
notice that  $ D(M-K)(L(i.)) \simeq Gr^{\cal W}_0 (\Omega^{*}
(Gr^{W^K}_{a_K(J,r)}  (L(i.))) $ is acyclic by the inductive
hypothesis in dimension $ < n$.
 This ends the proof of  the lemma.

\smallskip
\noindent {\em The relation between $ C^K_r L$ and $C^{KM}_r L$}\\
The following result describes $C^{KM}_r L$ as the fiber of an
intersection complex of the local system defined by $ C^K_r L$.

\smallskip
\noindent {\bf Proposition.} {\em Let $H = H^* (C_r^K L$) be
considered as a nilpotent orbit with indices $i \in M-K$, then we
have  a quasi-isomorphism: $C^{KM}_r L \cong {\cal W}_{-1
}\Omega^*(H)$.}

\smallskip
\noindent The proof is based on a decomposition as above: $(K_.^+
\times S(K) \times (M-K)_.^+ \times S(M-K)) \simeq M_.^+\times
S_K(M)$, with the correspondence $ (J,s.,J',s.')\rightarrow
((J,J'),(s'_.\cup K, s.)) $  then using the relations:

\noindent $a_{s_{\lambda}}((J,J'),r-1) = a_{s_{\lambda}}(J,r-1)$
when $ s_{\lambda}\subset K $

\noindent $a_{s'_{\lambda}}((J,J'),r-1) = a_{s'_{\lambda}}(J',-1)
+ a_K(J,r)$ and

\noindent $W^{s'_{\lambda}\cup K}_{a_{s'_{\lambda}\cup
K}((J,J'),r-1)} (Gr^{W^K}_{a_K(J,r)}) =
W^{s'_{\lambda}}_{a_{s'_{\lambda}}(J',-1)} (Gr^{W^K}_{a_K(J,r)}) $
since $W^{s'_{\lambda}\cup K}$ is relative to $W^K$\\
we find

 \noindent $C_r^{KM}L \colon=  s( \cap _{K\neq s_{\lambda}\in (s.'\cup
K, s.)}W^{s_{\lambda}}_{a_{s_{\lambda}}((J,J'),r-1)}
Gr^{W^K}_{a_K((J,J'),r)} L)_{((J,J'),(s'_.\cup K, s.)) \in
M_.^+\times S_KM} \simeq $

\noindent $s[ (\cap _{ s'_{\lambda}\in s.'}W^{K \cup
s'_{\lambda}}_{a_{s'_{\lambda}}(J', a_K(J, r-1)} (s [
Gr^{W^K}_{a_K(J,r)}(\cap _{K\supsetneqq s_{\lambda}\in s.}
W^{s_{\lambda}}_{a_{s_{\lambda}}(J,r-1)} L ]_{(J,s.)\in K_.^+
\times S(K)}) ]_{(J',s.')\in (M-K)_.^+ \times S(M-K)}$

\noindent $ \simeq s[ \cap _{ s'_{\lambda}\in s.'}
W^{s'_{\lambda}}_{a_{s'_{\lambda}}(J',-1)}  ( C_r^K
L)]_{(J',s.')\in (M-K)_.^+ \times S(M-K)} $

\noindent where we used  $\cap _{ s'_{\lambda}\in s.'}
W^{s'_{\lambda}}_{a_{s'_{\lambda}}(J',-1)}  ( C_r^K L) \simeq
C_r^K ( \cap _{ s'_{\lambda}\in s.'} W^{K \cup
s'_{\lambda}}_{a_{s'_{\lambda}}(J',-1)} L)$ is defined as above
for each subset  of $L$.

\noindent This formula is precisely the weight ${\cal W}_{-1}$ of
the terms of $C_r^{K}L$. It shows that $C_r^{KM}L$ is constructed
in two times, once as $C^K$ over $K_.^+ \times S(K)$ (that is a
space normal to $Y_K$) and once as a weight filtration over
$(M-K)_.^+ \times S(M-K)$
 (that is the space $Y_K$).
\section { Global construction of the weight filtration.}
In this section we construct
 a global bi-filtered combinatorial logarithmic complex and prove
a global decomposition of the graded weight into intermediate
extensions of polarised $VHS$ on the various intersections of
components of $Y$. We use a formula of the intersection complex
announced by Kashiwara and Kawai [26] that we prove since we have
no reference for its proof.

\smallskip
\n  Let $Y $ be a $NCD$  in $X$ with smooth irreducible components
  $Y_{i\in I}$ with indices
 in the set  $I$. The direct image of the complex local system $ {\cal L}^c$ is
 computed globally via the logarithmic complex with coefficients
  in Deligne's analytic extension
  ${\Omega}^*_X (Log Y) \otimes {\cal L}_{X}$. It is on a quasi-isomorphic
 constant combinatorial   complex with indices $s. \in S(I)$\\
\centerline {$ {\Omega}^* {\cal L} = s({\Omega}^*_{X_{s.}} (Log Y)
\otimes {\cal L}_{X_{s.}})_{s.\in S(I)} $}

\n that we can define the two filtrations $ {\cal W}$ and $F$.
  \subsection {Comparison with the
local definition} {\bf Lemma}. {\em Let $M \subset I $,  $y \in
Y^*_M$ and $L \simeq {\cal L}_X (y)$ the space of multivalued
sections of ${\cal L}$ at $y$, then the correspondence from $v \in
L$ to $\tilde v \in {\cal L}_{X,y}$ extends to a quasi-isomorphism
}

\smallskip \centerline {$ {\Omega}^* L \cong
 ({\Omega}^*{\cal L})_y$}

\n Proof. The quasi-isomorphism $\Omega (L,  N_j, j \in M ) \cong
({\Omega}^*_X (Log Y) \otimes {\cal L}_X)_y $ (2.1 (6)) is
compatible with the differentials when defined with indices $s.
\in S(I)$.

\smallskip
\n {\it The weight ${\cal W}$.}\\
1- For each $(s.) \in S(I)$ we deduce, from
 the weight filtration  by  sub-complexes $W_r(s.) L = s(W_r(J,s.)L)_{J\subset M}$
(3.2) of the locally defined
 DeRham complex $\Omega (L,N.)$,
 a corresponding global filtration
by sub-complexes\\
 ${\cal W}_r(X, {\cal L})(s.) $ in ${\Omega}^*_X (Log Y) \otimes {\cal L}_X$. \\
 In terms of a set of coordinates $z_i, i\in M $ defining $Y_M$
 in a neighbourhood of $y \in Y^*_M$, the fiber at $y$
  in $({\Omega}^*_X (Log Y) \otimes {\cal L}_X)_y$
  is defined as follows

  \smallskip
\n {\bf  Definition.} {\em ${\cal W}_r(X, {\cal L})_y(s.) $ is
generated as an ${\Omega}^*_{X,y}-$ sub-module
  by the germs of the sections $\wedge_{j \in J} \frac
{dz_j}{z_j} \otimes \tilde {v} $ for $v \in W_r(J,s.) L$.}

\smallskip
\n This formula is independent of the choice of coordinates, since
if we choose a different coordinate $z'_i = f z_i$ instead of
$z_i$ with $f$ invertible holomorphic at $y$, the difference $
\frac {dz'_i}{z'_i} - \frac {dz_i}{z_i} = \frac {df}{f}$ is
holomorphic at $y$, hence the difference of the sections
$\wedge_{j \in J} \frac {dz'_j}{z'_j}\otimes \tilde {v} -
\wedge_{j \in J} \frac {dz_j}{z_j} \otimes \tilde {v}$  is still a
section of the sub-complex
${\cal W}_r(X, {\cal L})(s.)$.\\
Moreover the restriction of the section is still defined in the
sub-complex near $y$, since $W_r(J,s.) L
\subset W_r(J-i,s.) L$ for all $i \in J$, then we have a quasi-isomorphism

\smallskip
\n {\bf  Lemma}. {\em We  have an induced quasi-isomorphism
$W_r(s.) L \cong  {\cal W}_r(X, {\cal L})_y(s.)$, functorial in
$(s.)$.}
 \be{Definition 1}{4.1} {\em (weight and
Hodge filtrations). The weight filtration is defined on the
combinatorial logarithmic complex with indices $s. \in S(I)$
$$ {\Omega}^* {\cal L} = s({\Omega}^*_{X_{s.}} (Log Y) \otimes
{\cal L}_{X_{s.}})_{s.\in S(I)}\leqno{(15)}$$
  \noindent as follows:
 $${\cal W}_r(X, {\cal L}) \colon= s( {\cal W}_r(X, {\cal L})(s.))_{s. \in
 S(I)} $$
The Hodge filtration $F$ is constant in $(s.)$ and deduced from
Schmid's extension to ${\cal L}_X$
$$F^p (s.) = 0 \rightarrow  F^p{\cal L}_X \ldots \rightarrow
 {\Omega}^i_{X_{s.}} (Log Y) \otimes
F^{p-i}{\cal L}_{X_{s.}}\rightarrow \ldots ;\quad F^p =
s(F^p(s.))_{s.\in S}\quad$$ }
 \n The fiber of Deligne's bundle  $ {\cal L}_X (y)$ at the point $y$
is identified with the space $L$ of multivalued sections $L$.

\smallskip
\noindent {\bf Definition 2} {\em (weight). With the same
notations, let $M \subset I, \vert M\vert  = p $ and  $y \in
Y^*_M$, then
 in terms of a set of $n$ coordinates $z_i, i\in [1,n]$ where we
identify $M$ with $[1,p]$, we write a section \\ $f =
(f^{s.})_{s.\in S}$ with $ f^{s.} = \Sigma_{J\subset M,J'\cap  M =
\emptyset}
 \frac {dz_J}{z_J}\wedge dz_{J'}  \otimes f^{s.}_{J,J'}$ s.t.$
f^{s.}_{J,J'}$ is not divisible by $y_j, j \in J$, then $f$ is in
$({\cal W}_r {\Omega}^*{\cal L})_y $ if and only if $
f^{s.}_{J,J'}(y)$ in $ {\cal L}_X (y) = L$ satisfy
$$ \forall J \subset M,\, f^{s.}_{J,J'}(y) \in
{\bigcap_{s_{\lambda}\in s.,s_{\lambda} \subset M}}
W_{a_{\lambda}(J,r)}^{s_{\lambda}} L $$}
 \noindent 3- {\it Definition  with residue.}
Let $W^Y$ denotes the weight along $Y$ on ${\Omega}^*_X (Log Y)$;
   choose an order on $I$ and an integer $m$,
 then  the residue morphism $Res_m$ of order $m$ is defined  on
 $W^Y_m ({\Omega}^p_X (Log Y))\otimes {\cal L}_{X}$
 with value in ${\Omega}^{p-m}_{Y^m} \otimes {\cal L}_{Y^m}$ on the disjoint union of
 intersections of $m$ components of $Y$ (the residue does not commute with differentials).
   For $M \subset I$ s.t.
 $Y_M \neq \emptyset$ and $\vert M \vert = m $ we deduce the
 residue
 $Res_M: W^Y_m ({\Omega}^p_X (Log Y))\otimes {\cal L}_{X}
 \to \Omega^{p-m}_{Y_M} \otimes {\cal L}_{Y_M}$ by composition of the residue
 morphism $Res_m$ with the obvious projection. At a point $y \in
 Y_M$, the morphism induced on the fiber with value in
 $({\Omega}^{p-m}_{Y^m} \otimes {\cal L}_{Y^m})(y)$ is denoted by $Res_{M,y}$.

\smallskip
 \noindent {\bf Definition 3}. {\em With the same notations, let $M \subset I, \vert
M\vert  = p $ and  $y \in Y^*_M$, then
  the fiber of the sub-analytic sheaf
$({\cal W}_r {\Omega}^*{\cal L})_y$  at $y$ is defined
successively on its intersection with\\ $ (W^Y_i
({\Omega}^p_{X_{s.}} (Log
Y)\otimes {\cal L}_{X_{s.}})_y$ by the following formula:\\
  a section $f \in (W^Y_i {\Omega}^p_X (Log Y) \otimes {\cal
L}_{X})_y$ is in $({\cal W}_r {\Omega}^*{\cal L})_y $ if and only
if setting ${\cal L}_{Y^m}(y) = L $}
 $$ \forall J, s_{\lambda} \subset M : \,\vert J
\vert = i,\, Res_{J,y} (f) \in {\Omega}^{p-m}(y) \otimes
W_{a_{\lambda}(J,r)}^{s_{\lambda}}  {\cal L}_{Y^m}(y) $$
 Remark.  By construction, for all
integers $r, \;{\cal W}_r / X-Y  = {\Omega}^*({\cal L})/ X-Y $, so
that ${\cal W}_r$ is exhaustive for $r$ big enough, and equal to
the extension by zero for $r$ small enough. It is a filtration by
sub-complexes of analytic sub-sheaves globally defined on $X$.

\smallskip
\noindent {\bf Proposition} ({\em Comparison with the local
definition). Let  $L = {\cal L}_X (y)$ denotes the space of
multi-valued sections of ${\cal L}$ at a point $y\in Y^*_M$, then
we have a bi-filtered quasi-isomorphism}

\smallskip \centerline {$ ({\Omega}^* L, {\cal W}, F) \cong
 ({\Omega}^*{\cal L}, {\cal W},
F)_y$}

\smallskip \n  At left we sum on $S(M)$ and at right on S(I).
The statement asserts that on $Y^*_M$, we still have a
quasi-isomorphism. For a fixed subset $M$ in $I$, the definition
of the weight filtration near a point in $Y^*_M$ involves only the
subsets $s_{\lambda} \subset M$, hence if we consider the
correspondence $s.' \in S(M)$ with the family $( M \cup s.'', s.')
\in S_M(I)$ where $s.'' \in S(I-M)$, the diagonal embedding of the
restrictions to $Y^*_M$,
 ${\Omega}^*_{X_{(s.')}} (Log Y) \otimes
{\cal L}_{X_{(s.')}},{\cal W}, F )_{\vert Y^*_M} $ into

\centerline {$s(({\Omega}^*_{X_{( M \cup (s.''), (s.'))}} (Log Y)
\otimes {\cal L}_{X_{( M \cup (s.''), ( s.'))}}, {\cal W},
F)_{\vert Y^*_M})_{ (s.'') \in S(I-M)}$}

\n is a bi-filtered quasi-isomorphism, hence the local study at
points of $Y^*_M$ of
 $(\Omega^{*} L,{\cal W}, F)$ reduces to
 $s({\Omega}^*_{X_{(s.')}} (Log Y) \otimes {\cal L}_{X_{(s.')}},
 {\cal W}, F )_{\vert Y^*_M})_{s.' \in S(M)}$.

\smallskip
\noindent {\bf The variation of Hodge structures }$({\cal L}^K_r,F)$ \\
\n Let $i_K:Y^{*}_K \rightarrow X $. Recall the definitions in the introduction:\\
 ${\cal L}^K  = i_K^* R^{\vert K \vert}{\bf j}_*
 {\cal L}$, \quad
  ${\cal L}'^K = i_K^* R^0 {\bf j}_* {\cal L}$ \\
   ${\cal L}^K_{X} =
 {\cal L}_{Y_K}/(\Sigma_{i \in K} {\cal N}_i {\cal L}_{Y_K}) $,
 \quad
${\cal L}'^K_{X} = \cap_{i \in K} (ker {\cal N}_i: {\cal L}_{Y_K}
\to {\cal L}_{Y_K}) $\\
and $ {\cal W}_{Y_K}^K = {\cal W}(\Sigma_{i\in K} {\cal N}_i)$
   for the filtration by sub-bundles defined on
${\cal L}_{Y_K}$ by $\Sigma_{i\in K} {\cal N}_i$.
 \be{Proposition.}{4.2}{\em i)
 ${\cal L}^K_{X} $  (resp. ${\cal L}'^K_{X}$)
 induces a flat bundle on $Y^{*}_K $, with flat sections isomorphic to
  the  local system
 ${\cal L}^K$ (resp. ${\cal L}'^K$); precisely  they are respectively Deligne's
extension of the corresponding complex local system.

\smallskip
\n ii) The filtration $ {\cal W}_{Y_K}^K $ induces a filtration by
flat sub-bundles of ${\cal L}^K_{X} $ (resp ${\cal L}'^K_{X} $) on
$Y^{*}_K$, hence induces a filtration by complex sub-local systems
$ {\cal W}^K $ of ${\cal L}^K$ ( resp. ${\cal L}'^K $).

\smallskip
\n iii) The filtration $ {\cal W}^K $ is rationally defined on the
rational local system $ {\cal L}^K $ ( resp. ${\cal L}'^K $).

\smallskip
\n iv)  Let ${\cal L}_r^K \colon= Gr^{{\cal W}^K}_{r - \vert K
\vert }{\cal L}^K$ and ${\cal L}_{X,r}^K \colon= Gr^{{\cal
W}^K_{Y_K}}_{ r -
\vert K \vert }{\cal L}^K_{X}$ for $r>0$, resp.\\
 ${\cal L}_r^K \colon= Gr^{{\cal W}^K}_{r + \vert K \vert }{\cal L}'^K$
 and ${\cal L}_{X,r}^K \colon=
Gr^{{\cal W}^K_{Y_K}}_{ r + \vert K \vert }{\cal L}'^K_{X}$ for
$r< 0$,\\
 then the  system defined on $Y_K^*$ by  $({\cal L}_r^K,
 ({\cal L}_{X,r}^K)_{\vert Y_K} , F)$
 where $F$ is the Hodge filtration induced from ${\cal L}_X$,
  is a polarised variation of Hodge structures $VHS$ of
 weights
$r - \vert K \vert + m$ for $r > 0$ and $r + \vert K \vert + m$
for $ r < 0 $.}

\smallskip
\noindent Proof. We deduce from the comparison propositions  with
local definitions at each point $y \in Y^*_K$ the following
complexes quasi-isomorphic to ${\bf j}_*{\cal L}_y$

\smallskip \centerline {$ ({\Omega}^*_X (Log Y) \otimes
{\cal L}_X)_y \cong
 \Omega (L,  N_j, j \in K )\cong s(L(J), N.)_{J \subset
K}, \;({\Omega}^* L, {\cal W}, F) \cong
 ({\Omega}^*{\cal L}, {\cal W}, F)_y$}

\n  hence near each point $y \in Y_K^*$, ${\cal L}^K$  and ${\cal
L}'^K$ are (locally) constant.  The local system ${\cal L}^K$
(resp. ${\cal L}'^K$) is defined by the flat sections of the
bundles ${\cal L}^K_X$ (resp. ${\cal L}'^K_X $ ) whose connection
has logarithmic singularity since it is induced by the connection
on ${\cal L}_X$ which proves (i). \\
The same argument apply to the filtration $ {\cal W}_{Y_K}^K $,
which proves  ii).

\smallskip
\noindent iii) Let us denote by ${\cal L}^{K,rat}$ the rational
local system underlying the complex  ${\cal L}^{K}$. The
intersection $ {\cal W}^K \cap {\cal L}^{K,rat} $   defines a
rational filtration underlying the complex one. This can be
checked locally as the graded vector space $ Gr^{W^K}_r  L $ has a
rational structure at each point $y$.

\smallskip
\noindent iv) The sheaf  ${\cal L}^K _r$ is locally constant and
isomorphic to the cohomology  of the complex $C^K_{r} L$ for $L =
{\cal L}_X (y)$ (Prop. 3.6) which shows that the local system
${\cal L}^K_r$ is defined by the flat sections of the bundle
${\cal L}_{X,r}^K$, then (iv) follows.

\smallskip
\noindent {\it Remark}. Given the $VHS$ ${\cal L}^K_{r}$, we can
construct a corresponding complex $({\Omega}^*{\cal L}^K_{r} ,
{\cal W}, F)$, then for each point $y \in Y_M^*$, $K \subset M$,
we have a quasi-isomorphism $ C^{KM}_r  L \cong {\cal
W}_{-1}({\cal L}^K_{r} [-\vert K \vert])_y $ for $r > 0$ and
  $ C^{KM}_r L \cong {\cal W}_{-1}({\cal L}^K_{r} [1-\vert K \vert])_y$ for $r <
0$ (recall $ C^{KM}_0  L \cong 0$). That is the cohomology of the
various complexes $ C^{KM}_r  L $ is globally defined on $Y_M^*$.
In fact we will see in the next result it is the restriction of
the intermediate extension of ${\cal L}^K_{r}$ to $Y_M^*$.

\smallskip
\noindent it will follow from the next proof that we could define
${\cal L}^K_r$ as \\
 ${\cal L}^K_r \colon= i_K^* {\cal H}^{\vert K
\vert}(Gr^{\cal W}_{r}{\Omega}^* {\cal L}), r > 0  \,\, (resp.\,
\,{\cal L}'^K_r \colon= i_K^* {\cal H}^{\vert K \vert -
1}(Gr^{\cal W}_{r}{\Omega}^* {\cal L} ), r < 0)$
  \be {Theorem.}{4.3}{\em Let ${\cal L}$ be a  local system  with
locally unipotent monodromy, underlying a variation of polarised
Hodge structures of weight $m$ on $X-Y$ of $dim \; n$ and let
$j^K:Y^{*}_K \rightarrow Y_K $, $i_K:Y^{*}_K \rightarrow X $, then
the bi-filtered complex

\smallskip
\centerline {$ ({\Omega}^*{\cal L}, {\cal W}, F)$}

\smallskip
 \noindent is filtered
quasi-isomorphic to $ ({\Omega}^*_X (Log Y) \otimes {\cal L}_X, F)
$.\\
i) The restriction to ${Y_K^*}$, $i_K^* {\cal H}^{\vert K
\vert}(Gr^{\cal W}_{r}{\Omega}^* {\cal L}) $ for $r > 0$ (resp.
$i_K^* {\cal H}^{\vert K \vert - 1}(Gr^{\cal W}_{r}{\Omega}^*
{\cal L} )$ for $r < 0$) is a complex local system isomorphic to $
{\cal L}^K_r$, moreover the following decomposition property into
 intermediate extensions is satisfied

\smallskip
\noindent  \centerline {$(Gr^{\cal W}_{r}{\Omega}^* {\cal L}, F)
\cong \oplus_{K\subset I} {\bf  j}_{!*}^K {\cal L}^K_{r} [ - \vert
K \vert], F[-\vert K \vert])$, for $ r > 0$ }

\smallskip
\centerline {$(Gr^{{\cal W}}_{0}{\Omega}^* {\cal L}, F) \cong 0$}

\smallskip
\centerline {$(Gr^{{\cal W}}_{r}{\Omega}^* {\cal L}, F) \cong
\oplus_{K \subset I} j_{!*}^K {\cal L}^K_{r} [1 - \vert K \vert ],
F) $, for $ r < 0$}

\smallskip
\noindent ii) (Kashiwara and Kawai's formula [26]) The sub-complex
$ {\cal W}_{-1}{\Omega}^*({\cal L}[2n])$ is quasi-isomorphic to
the intermediate extension ${\bf j}_{!*}{\cal L}[2n]$ of
 ${\cal L}[2n]$ and $ ({\cal W}_{-1}{\Omega}^*({\cal L}[2n]),F)$ is a Hodge
complex of weight $2n + m$.}

\smallskip
 \noindent
 Proof. i) The decomposition of $(Gr^{{\cal W}}_{r}
{\Omega}^* {\cal L} ,F) $ reduces near a point $y \in Y^*_M $ to
the local decomposition for the nilpotent orbit $L$ defined at the
point $y$ by the local system
 $ Gr^{\cal W}_{r}{\Omega}^* L \cong
\oplus_{K \subset M}C_r^{KM}L $ . The global decomposition that
follows is $Gr^{\cal W}_{r}{\Omega}^* {\cal L} \cong
\oplus_{K\subset I}{\cal W}_{-1}{\Omega}^* ( {\cal L}^K_{r})$ that
is the complex ${\Omega}^*$ considered for the polarised local
system ${\cal L}^K_{r}$, as it follows from the local formula
$C^{KM}_rL \cong {\cal W}_{-1} L^K_r)$ where $ L^K_r \simeq
H^*(C^K_rL)$ which has been checked. The fact that $C^{KM}_rL $ is
precisely the fiber of $j_{!*}^K {\cal L}^K_{r} [-\vert K\vert ]$
for $r
>0$ (resp. $j_{!*}^K {\cal L}^K_{r} [1 -\vert  K\vert ]$ for $r <
0$ will follow from (ii) by induction on the dimension.\\
The count of weight and the shift in $F$ take into account for $r
> 0$ the residue in the isomorphism with $L$ that shifts $W$ and $
F$ but also the shift in degrees, while for $ r < 0 $ there is no
residue (since the cohomology is in degree $0$ of the logarithmic
complex with index $s. = K \in S(K)$ but only a shift in degrees
$\vert K \vert - 1$ in the combinatorial complex, the rule being
as follows:

\n {\em Let $( K, W, F)$ be a mixed Hodge complex then for all $m,
h \in \Z$, $( K', W', F') = ( K [m], W[m-2h], F[h])$ is also a
mixed Hodge complex.}\\
 The same proof apply for $r =0$.

\smallskip
{\it   ii)} The proof is based on the decomposition of $S_KM$ as a
product in (3.3) and follows by induction on the dimension $n$,
from the local decomposition of the graded parts of
the weight filtration above in i). \\
 The proof  is true in
dimension $1$ and if we suppose the result true in dimension
strictly less than $n$, we can apply the result for $Gr^{\cal
W}_{r}{\Omega}^* {\cal L} $ that is for local systems defined on
open subsets of the closed sets $Y_K$, namely the local system
${\cal L}^K_r [-\vert  K\vert ]$ for $r > 0$ (resp.${\cal L}^K_r
[1 -\vert  K\vert ]$ for $r < 0$) whose fiber at each point $y \in
Y^*_K$ is quasi-isomorphic to $C^K_rL$. Let $ j^K:Y^{*}_K
\rightarrow Y_K $ be  the open embedding in $Y_K$ and consider the
associated DeRham complex ${\Omega}^*({\cal L}^K_r)$ on $Y_K$
whose  weight filtration will be denoted locally near a point in
$Y^*_M$ by ${\cal W}^{M-K}$ for $K \subset M$; then by the
inductive hypothesis we have at
 the point $y $:  ${\cal W}^{M-K}_{-1} {\Omega}^* L^K_r \simeq {\cal W}^{M-K}_0
{\Omega}^* L^K_r  $ is also quasi-isomorphic to the fiber of the
intermediate extension of ${\cal L}$, that is
 $$ \forall r > 0, C_r^{KM}(L)\simeq (j_{!*}^K {\cal L}^K_{r}
  [-\vert  K\vert ])_y \simeq {\cal W}^{M-K}_{-1 }{\Omega}^* L^K_r $$
and similarly for $r < 0$.

In order to check  the result for $ {\cal W}_0{\Omega}^*({\cal
L}[2n])$, we use the following criteria characterising
intermediate extension [17] where the degree shift is by $2n$:

\noindent Consider the stratification defined by $Y$ on $X$ and
the middle perversity $p(2k) = k - 1 $ associated to the closed
subset $Y^{2k} = \cup_{\vert K\vert  = k} Y_K$ of real codimension
$2k$. We let $Y^{2k-1} = Y^{2k} $ and $p(2k-1) = k -1 $. For any
complex of sheaves $S.$ on $X$ which is constructible with respect
to the stratification, let $S.^{2k} = S.^{2k-1}  = S.\vert
X-Y^{2k}$ and consider the four properties:

a) Normalization: $S.\vert  X-Y^2 \cong {\cal L}[2n]$

b) Lower bound: ${\bf H}^i (S.) = 0 $ for all $i < -2n$

c) vanishing condition: ${\bf H}^m ( S.^{2(k+1)})= {\bf H}^m (
S.^{2k +1})= 0  $ for all $m > k -2n$

d) dual condition: ${\bf H}^m (j_{2k}^! S.^{2(k+1)})= 0  $ for all
$k \geq 1$ and all $m > k -2n$ where  \\
$j_{2k} : (Y^{2k} -
Y^{2(k+1)}) \rightarrow (X - Y^{2(k+1)})$ is  the closed
embedding,

 \noindent then we can conclude that $S.$ is the intermediate extension of
 ${\cal L}[2n]$.

\noindent In order to prove the result in dim. $n$ we check the
above four properties for $ {\cal W}_0{\Omega}^*({\cal L}[2n])$.
The first two are clear and  we use the exact sequences $$ 0
\rightarrow {\cal W}_{r-1}{\Omega}^*({\cal L}[2n])\rightarrow
{\cal W}_r {\Omega}^*({\cal L}[2n])\rightarrow Gr^{\cal W}_r
{\Omega}^*({\cal L}[2n]) \rightarrow 0 $$ to prove d)(resp. c)) by
descending (resp. ascending )indices from ${\cal W}_r $ to ${\cal
W}_{r-1}$ for $r \geq 0$ (resp. $r-1 $ to $r$ for $r < 0 $ )
applying at each step the inductive hypothesis to $Gr^{\cal W}_r
$.

Proof of d). For $r$ big enough  $ {\cal W}_r {\Omega}^*({\cal
L}[2n])$ coincides with the whole complex $ {\bf j}_* {\cal
L}[2n]$, then the dual condition is true for $r$ big enough. Now
to check d) for $Gr^{\cal W}_r {\Omega}^*({\cal L}[2n])$,  we
apply d) to a component with support $Y_{K'} $  with $\vert
K'\vert \,= \,k'$. We choose  $k
> k'$ and consider   $j'_{2k} : (Y^{2k}\cap
Y_{K'} - Y^{2(k+1)}\cap Y_{K'}) \rightarrow (Y_{K'}  -
Y^{2(k+1)}\cap Y_{K'} )$  ( notice that $Y^{2k}\cap Y_{K'} = (Y
\cap Y_{K'})^{2(k-k')}$ is of codim. $2(k-k')$ in $Y_{K'}$),
 then for $S.' $ equal to the intermediate extension of
 ${\cal L}^{K'}_r [2n - 2k']$ on $Y_{K'}$ we have on $Y_{K'}$ the property
 ${\bf H}^m ({j'}^!_{2k} S.'^{2(k-k')+1)})= 0
$ for all $ (k - k') \geq 1$ and all $m > k - k' - 2(n-k') = k+k'
- 2n$ which gives for $ S.' [k'] $ on X: ${\bf H}^m (j_{2k}^!
S.'^{2(k+1)}[k'] )= 0  $ for all $k > k' $ and all $m > k  - 2n$,
hence d) is true.\\
 If $k = k'$, then $Y^{2k}\cap Y_{K'}  =
Y_{K'}$ and we have a local system in degree $k'-2n $ on $Y_{K'} -
Y^{2(k+1)}\cap Y_{K'} $ hence d) is still true, and for $k < k' $,
the support   $Y_{K'}  - Y^{2(k+1)}\cap Y_{K'} $ of $S.'$ is
empty. From the decomposition theorem and the induction, this
argument apply to $Gr^{\cal W}_r $ and hence apply by induction on
$r \geq 0$ to $ {\cal W}_0 $ and also to ${\cal W}_{-1} $.

Proof of c). Dually, the vanishing condition is true for $r$ small
enough since then $  {\cal W}_r $ coincides with the  extension by
zero of ${\cal L}[2n]$ on $X-Y$.

\noindent Now we use the filtration for $r < 0 $, for $S.' $ equal
to the intermediate extension of ${\cal L}^{K'}_r [2n - 2k']$ on $
Y^{K'}$ we have for $k > k'$: ${\bf H}^m (S.'^{2(k-k')+1})= 0  $
for all $m > k + k' - 2n$, which gives for $ S.' [k'+1],\; (r < 0)
$ on X: ${\bf H}^m ( S.^{2(k+1)})= {\bf H}^m ( S.^{2k +1})= 0  $
for all $m > k-1 -2n$. If $k = k'$, then $ S.' [k'+1] $ is a local
system in degree  $-2n + k-1 $ on $Y_{K'} - Y^{k+1}$ and for $k <
k' $, $Y_{K'} - Y^{k+1}$ is empty.

\smallskip
\noindent { \bf Corollary}. {\em The weight filtration $ {\cal W}$
of ${\Omega}^* {\cal L}$ is defined over $\Q$}.

\smallskip
\noindent The proof is based on the following lemma applied to
${\bf j}_* {\cal L}^c \simeq {\Omega}^* {\cal L}$ with its
filtration ${\cal W}$.

\smallskip
\noindent {\bf Lemma}. {\em Let $K$ be a $\Q-$perverse sheaf such
that $K^c = K \otimes \C$ is filtered by a finite filtration $W^c$
of complex perverse sub-sheaves s.t.$Gr^{W^c}_r K^c$ is rationally
defined and  the rational filtration $ W^r = W^c \cap K$ induces
the rational structure on $Gr^{W^c}_r K^c$, then $ W^r$ is a
rational filtration by perverse sub-sheaves of $K$ such that $W^r
\otimes \C \simeq W^c$.}

\smallskip
\noindent The proof is similar to the case of local systems and is
by induction on the weight $i$ since by hypothesis it applies to
the lowest weight. Considering the extension
 $0 \to W^c_i K^c \to W^c_{i+1} K^c {\stackrel \pi \to} Gr^{W^c}_{i+1}K^c
\to 0$, then  $ W^r_{i+1} K \otimes \C \simeq W^c_{i+1} K$ follows
from the hypothesis $(Gr^{W^r}_{i+1}K) \otimes \C \simeq
(Gr^{W^c}_{i+1}K)$ and the inductive isomorphism for  $ W^r_{i}
K$.

\smallskip
\noindent { \bf Corollary}. {\em
 If  $X$ is proper and if we forget the negative weights in   the filtration
${\cal W}$ that is we consider  ${\cal W}''$ with ${\cal W}_i'' =
{\cal W}_i $ for $i \geq 0$ and ${\cal W}_{-1}'' = 0$, then the
bi-filtered complex
$$ ({\Omega}^*({\cal L}), {\cal W}''[m], F)$$
is a  mixed Hodge complex .}

\section{The complex of nearby cycles  ${\Psi}_f (\cal L)$.}
Let $f \colon X \rightarrow D$ and suppose $Y = f^{-1}(0)$ a
$NCD$, the complex of sheaves $\Psi_f {\cal L}$ of nearby
co-cycles on $Y$ has been introduced in [11]; its cohomology fiber
$H^i ( ({\Psi}_f {\cal L)}_y) \simeq H^i (F_y , {\cal L})$ at a
point $y$ in $Y$ is isomorphic to the cohomology of the Milnor
fiber $F_y$  at $y$. The monodromy  induces an action ${\cal T}$
on the complex itself. If $dim \, X = n$, $\Psi_f {\cal L} [n-1]$
is perverse on $Y$. Since the local system ${\cal L} $ is defined
over $\Q$, the monodromy decomposes in the abelian category of
$\Q-$perverse sheaves as the product ${\cal T} = {\cal T}^s {\cal
T}^u$ of  simple and unipotent endomorphisms. \noindent  Let
${\cal N} = Log{\cal T}^u$, then Deligne's filtration ${\cal W}
({\cal N} ) $ is defined over $\Q$.
 The aim of this section is to describe the structure of a mixed
  Hodge complex ($MHC$)
on ${\Psi}_f {\cal L}$ with weight filtration ${\cal W} ({\cal N}
) $.
 This problem is closely related to the weight
filtration in the open case since we have the following {\it
relation between  ${\Psi}^u_f {\cal L}$,  the direct image $ {\bf
j}_* {\cal L}$ and ${\bf j}_{!*} {\cal L}$ } as  explained in [1]
and [2]
$$ i_Y^*({\bf j}_* {\cal L}[n] / {\bf j}_{!*} {\cal L}[n]) \quad \simeq \quad Coker
({\cal N} \colon {\Psi}_f  ^u {\cal L} [n-1]\rightarrow {\Psi}_f^u
{\cal L}[n-1])  \leqno{} $$ \noindent The  filtration $W({\cal
N})$ on ${\Psi}_f^u ({\cal L})$  induces a filtration ${\cal W}$
on $Coker ({\cal N} / {\Psi}_f^u {\cal L})$, hence on ${\bf j}_*
{\cal L} / {\bf j}_{!*}
{\cal L}$.\\
{\it  The induced  filtration on ${\bf j}_* {\cal L} / {\bf
j}_{!*}  {\cal L}$  is  independent of the choice of $f$.} For a
rigorous proof one should use the result of Verdier [34].
 A path in the space of
functions between two local equations $f$ and $f'$ of $Y$ gives
rise to
 an  isomorphism between  ${\Psi}_f {\cal L}$ and
${\Psi}_{f'} {\cal L}$ ; it is only modulo $Coker {\cal N}$ that
this isomorphism is independent of the path. We do check here that
the  weight filtration $ {\cal W}$ on ${\Omega}^* {\cal L}$ is
induced locally by $W( {\cal N})$.
\subsection{ The weight filtration on the nearby co-cycles ${\Psi}_f {\cal
L}$}

When we consider the coefficients in the complex local system
${\cal L} \otimes \C$ ( denoted also ${\cal L}$), the method to
compute ${\Psi}_f$ as explained in [11] uses the restriction
$i^*_Y \hbox {\bf j}_* {\cal L}$ of the higher direct image  of
${\cal L}$ to $Y$ and the cup-product $H^i(X^*, {\cal L}) \otimes
H^1(X^*, \Q ) \; {\buildrel {\smile  \eta } \over \longrightarrow}
\;H^{i+1}(X^*, {\cal L})$ by the inverse image $\eta = f^* c \in
H^1 (X^*, \Q )$ of a  generator $c$ of the
 cohomology $H^1(D^*, \Q )$.
 We construct effectively, using Deligne's bundle extension, a bi-filtered complex on
 which $\eta$ is defined as a
morphism (of  degree 1), $\eta \colon \; i^*_Y ({\Omega}^*_X (Log
Y) \otimes {\cal L}_X) \to i^*_Y ({\Omega}^*_X (Log Y) \otimes
{\cal L}_X) [1 ]$ satisfying ${\eta}^2 = 0$ so to get a double
complex whose simple associated complex is quasi-isomorphic to
${\Psi}^u_f ({\cal L})$.
 \be{The global weighted complex
$({\Psi}^u_f ({\cal L}), {\cal W}, F)$} {5.1} \

Let $t$ denotes a coordinate on the disc $D$ and $\eta \colon= f^*
(\frac{dt}{t})$, then $\wedge \eta$ defines a morphism of degree
one on $i_Y^{*}{\Omega}^*_X (Log Y) \otimes {\cal L}_X$. We
consider the simple complex
$$ \Psi^u_f {\cal L}_X \colon=  s(i_Y^{*}({\Omega}^*_X (Log Y)
 \otimes {\cal L}_X )[p],\eta)_{p \leq 0}\leqno{(16)} $$
defined by the double logarithmic complex ( $\oplus_{p \leq
0}i_Y^{*} ({\Omega}^i_X (Log Y) \otimes {\cal L}_X)$ is in degree
$i$). To define as previously a constant combinatorial resolution
of $\Psi^u_f {\cal L}_X $, we put $ {\Psi}^{u}_f {\cal L}_X(s.) =
\Psi^u_f {\cal L}_X$ for each $s. \in S(I)$ and let
$$ {\Psi}^{u}_f {\cal L} \colon= s({\Psi}^{u}_f {\cal L}_X(s.))_{s. \in S(I)}
\simeq s(i_Y^{*}{\Omega}^* {\cal L} [p],
 \eta)_{p \leq 0}\leqno{(17)}$$
 which can be viewed also as $s({\Omega}^* {\cal L} [p],
 \eta)_{p \leq 0}$;
then we define  the weight filtration and the Hodge filtration by
$$ {\cal W}_r ({\Psi}^{u}_f {\cal L}) =  s (i_Y^{*}{\cal W}_{r+2p-1}
 {\Omega}^* {\cal L} [p], \eta)_{p \leq 0}, \;
F^r ({\Psi}^{u}_f {\cal L}) =  s ( i_Y^{*} F^{r+p}
 {\Omega}^* {\cal L}[p], \eta)_{p \leq 0}\leqno{(18)}$$
The logarithm of the monodromy ${\cal N}$ is defined on this
complex and we want to show that  the filtration ${\cal W}$ above
coincides with  $W({\cal N})$.
 \be { Theorem.}{5.2} {\em Suppose
that $\cal L$ underlies a unipotent variation of polarised Hodge
structures of weight $m$, then $W({\cal N})= {\cal W}$.}

\smallskip
\noindent With this result we can conclude that the weight
filtration in the open
 case is induced locally by the weight filtration defined by the monodromy
  on the nearby co-cycles.\\
 The proof of this theorem is based on the results in
  the local case, that is for the nilpotent
 orbit defined by ${\cal L}$ at a
point  $y \in Y^*_M$, $M \subset I$.
  \be{Local description of the
weight and Hodge filtrations }{5.2}

\n Near a point $y \in Y^*_M, M \subset I$, we can find
coordinates $z_i$ for $i\in M$  defining $Y^*_M$ locally
 and non zero integers $n_i$ s.t. $f = {n \atop {\buildrel \prod \over {i=1}}}
z^{n_i}_i$ where  we do suppose $i \in [1,n]$, where $\vert M
\vert = n $ ( $n$ is less or equal to the $dim$ of $X$),
 then  in  DeRham cohomology  $\eta =
f^*(\frac{dt}{t}) = {n \atop {\buildrel \Sigma \over {i=1}} }  n_i \frac {dz_i} {z_i}$.\\
Thus $\eta$ defines a morphism of degree one on the DeRham complex
$ \Omega (L ,  N_i )_{i\in [1,n]}$ satisfying  ${\eta}^2 =0$. We
define
$${\Psi}^{0}(L) = s(\Omega (L ,  N_i)_{i\in [1,n]})[p],\eta )_{p \leq 0}$$
as the simple complex defined by the double complex
for $p \leq 0$. \\
{\em Remark}: In order to take into account the action of $ N =
Log T^u$ we may write $L[N^p]$ for $ L[p]$ and $ L [N^{-1}]$ for
the direct sum over $p \leq 0$, so that the action of $ N $ is
just multiplication by $N$, then
 $${\Psi}^{0}(L)  \simeq  \Omega (L[N^{-1} ],  N_i - n_i N)_{i\in [1,n]}. \leqno{} $$
is the Koszul complex on $ N_i - n_i N$ acting on $L[N^{-1}]$.\\
{\it The complex ${\Psi}^{0}_M L $}. To describe the weight in
terms of the filtrations $ ({\Omega}^* L, {\cal W}, F)$
 associated to $L$, we need to use the constant complex with index $s. \in S(M)$,
  ${\Psi}^{0} L(s.) = {\Psi}^{0} L $ and  introduce the complex
$$ {\Psi}^{0}_M L \colon= s({\Psi}^{0} L (s.))_{s. \in S(M)} \leqno{(19)} $$
 which can be viewed also as $ s({\Omega}^* L [p], \eta)_{p \leq 0}$,
 then we define on it the {\it weight and Hodge filtrations}
$$ {\cal W}_r ({\Psi}^{0}_M L) =  s ({\cal W}_{r+2p - 1}
 {\Omega}^* L [p], \eta)_{p \leq 0}, \; \;
F^r ({\Psi}^{0}_M L) =  s ( F^{r+p} {\Omega}^* L [p], \eta)_{p
\leq 0}.  \leqno{(20)}$$ \noindent {\it Monodromy.}

 The logarithm ${\cal N}$ of the monodromy is defined by an endomorphism
$\nu$ of the  complex $ {\Psi}^0_M L $, given by the formula

$$\forall a. = \Sigma_{p\leq 0} a_p \in  {\Psi}^0_M L: \quad (\nu
(a.))_p = a_{p-1} $$ such that $\nu ( {\cal W}_r ) \subset {\cal
W}_{r-2}$ and $\nu ( F^r ) \subset F^{r-1}$.

\smallskip
\noindent {\bf Lemma 1}. {\em The local quasi-isomorphism on the
stalk at a point $y \in Y^*_M$ of the logarithmic complex with
coefficients in Deligne's extension $(\Omega^*_X (Log Y) \otimes
{\cal L}_X)_y $ extends to a quasi-isomorphism from $(\Psi^u_f
{\cal L}_X)_y $ (16) to    ${\Psi}^0 L $
 (resp. from $({\Psi}_f^u {\cal L})_y $ (17)
 to $ {\Psi}^0_M L $) respecting the weight
and Hodge filtrations }
$$(({\Psi}^u_f {\cal L}_X)_y, W,F) \; \simeq ({\Psi}^0  L, W,F)
\quad , \quad (({\Psi}^u_f {\cal L})_y,{\cal W},F) \simeq
({\Psi}^0_M L,{\cal W}, F) $$ This lemma is the needed link
between the global and local cases.

\smallskip
\noindent {\bf Lemma 2}. {\em We have a triangle in the derived
category  represented by the exact sequence\\ $0 \rightarrow
i_Y^{*}({\Omega}^*_X (Log Y)
 \otimes {\cal L}_X  \rightarrow \Psi^u_f {\cal L}_X  {\stackrel
\nu \rightarrow}
 \Psi^u_f {\cal L}_X  \rightarrow 0$}
\subsection{ Main local results }
\noindent {\it Proof of  the theorem.} The proof can be reduced to
the local case at a point $y \in Y^*_M$. We introduce first the following complexes.\\
 The morphism $\eta $ induces a morphism denoted also by $\eta : C^{KM}_r
L \rightarrow C^{KM}_{r+2} L [1]$ so that we can define a double
complex and the associated simple complex
$$ \Psi^{KM}_r L = s( C^{KM}_{r+2p-1} L [p], \eta )_{p \leq 0}, \;
\Psi^{K}_r L \colon= \Psi^{KK}_r L \leqno{(21)}$$
 \n i) {\it Decomposition of $ Gr^{\cal W}_r ({\Psi}^{0}_M L) $.}
  {\em There exist  natural injections
of $\Psi^{KM}_r L $ into $ Gr^{\cal W}_r {\Psi}^0_M L $ and a
decomposition}
$$  Gr^{\cal W}_r \Psi^{0}_M L  \simeq \oplus_{K \subset M}
\Psi^{KM}_r L \leqno{(22)}$$ Proof:  The result follows from the
decomposition of $Gr^{\cal W}_r \Omega^* L $ in the previous  open
case, applied to the spectral sequence with respect to $p$ in the
double complex above.

\smallskip
\noindent ii) We introduce  now the  the complex $A^{K}_i $ and
prove
 \be { Basic lemma. }{} {\em For all $i \geq 1$, the complex
$A^{K}_i \colon= s[C^K_{i+2p-1} L [p], \eta]_{1-i \leq p \leq 0}
\cong 0$ is acyclic.}

\smallskip
\n Proof: We view $A^{K}_i$ as a double complex where $\eta$ of
degree $1$ is a differential of the direct sum of complexes
$C^K_{i+2p-1} L $ without shift in  degrees  :

\smallskip
\noindent \centerline { $A^{K}_i \colon= C^K_{-i+1}{ \stackrel
\eta \to } \cdots { \stackrel \eta \to } C^K_{i+2p-1} L \cdots
 {\stackrel \eta \to } C^K_{i-1} L, \quad  A^{K}_i (J) =
\oplus_{1-i \leq p \leq 0} C^K_{i+2p-1} L (J), J \subset K $}

 \n  We may
filter $A^{K}_i$ by sub-complexes $U_i$. One way is to take $U_0 =
C^K_{i-1} L$, $U_1 = s[C^K_{i+2p-1} L [p], \eta]_{1-i < p \leq 0}$
and $U_2 = A^{K}_i$, that is we write $ A^{K}_i$ as:

\smallskip
\noindent \centerline { $s [ C^K_{-(i-1)} L [-2], s[C^K_{i+2p-1} L
[p], \eta]_{0 < p <i-1 }[-1] , C^K_{i-1} L, \eta] $}

\noindent  so to use by induction $ Gr^U_1 =  A^{K}_{i-2} \cong
0$. \\
 However we will use the technique of the spectral sequence defined by the
  increasing filtration   \\
$U_r =  s[C^K_{i+2p-1} L [p], \eta]_{- r  \leq  p \leq 0}$ for
$0 \leq r \leq i-1 $ \\
 with Deligne's notations

 \smallskip
\noindent \centerline {$E^{a,b}_1 = H^{a+b}(Gr^U_{-a}) =
H^{a+b}(C^K_{i+2a-1} L) $}

\n Since $ H^u(  C^K_{j} L) = 0$ if $j <0$ and $ u \neq \vert  K
\vert-1$
 or $j > 0$ and $ u \neq \vert  K \vert$,
  we get\\
 1) if $2 a > 1-i$,
$E^{a,\vert  K \vert-1-a}_1= H^{\vert  K \vert-1}(C^K_{i+2a-1} L)=
Gr^{W^K}_{ \vert  K \vert+i+2a-1}[(\cap_{i \in K} ker N_i: L
\rightarrow
L ]  $ \\  and $0$ otherwise\\
2) if $2a < 1-i$, $E^{a,\vert  K \vert-a}_1 = H^{\vert  K
\vert}(C^K_{i+2a-1} L)= Gr^{W^K}_{- \vert  K \vert+i+2a-1}[ L
/(\Sigma_{i \in K}
N_i L)]$\\ and $0$ otherwise\\
3) if $2a = i-1$, $E^{a,b}_1= 0$.

Starting at the level $1$, the term  $E^{-(i-1)-d,\vert  K
\vert-1+(i-1+d)}_{1}$ remains unchanged until we reach the level
$r = i-1+2d $ where the only non trivial differential appears and
we want  to show it is an isomorphism. The proof is based
 on the following  study of this differential.

\n The cohomology space of $C^K_r L$ for  will be identified with
the following polarised subspace of $Gr^{{W^K}}_{r-\vert  K
\vert}L$. For each $r>0$, and $(m_1, \ldots, m_n) \in T(r)$ ( that
is $m_i > 1, \Sigma m_i = r + \vert  K
\vert$), let\\
$P(m.)L = \cap_{i\in [1,n] }(ker N_i^{m_i-1}:Gr^{W^n}_{m_n-2 }
\cdots Gr^{W^1}_{m_1-2 }L \to Gr^{W^n}_{-m_n } \cdots
Gr^{W^1}_{-m_1 }L) \subset Gr^{{W^K}}_{r - \vert  K
\vert} L$,\\be the  primitive polarised subset, then we have the isomorphism \\
$P(m.)L \simeq Gr^{W^n}_{m_n-2 } \cdots Gr^{W^1}_{m_1-2 }(L
/(\Sigma_i N_i L))$.

 \smallskip
\noindent {\bf Sub-lemma.} {\em Let $N_K = \Sigma_{i\in K} n_i
N_i$ and for each $ (m.) \in T(i+2d-1)$, let $ P(m.) L $ denotes
the primitive sub-$HS$ as above. The differential at the level $r
= i-1+2d $ of the spectral sequence

 \smallskip
\noindent \centerline {$E^{-(i-1)-d,\vert  K \vert-2+i+d}_{r} \to
E^{d,\vert K \vert -d}_{r}, \quad r = i-1+2d $}

\n is given by the inverse of the isomorphism up to a constant \\
$ Gr^{W^K}_{- \vert  K \vert+r}[ L /(\Sigma_{i \in K} N_i L)]
\simeq  \oplus_{m. \in T(r)} P(m.)L
 {\stackrel { N_K^{r - \vert  K \vert} }\longrightarrow }\oplus_{m. \in
T'(r)} P(m.)L \simeq Gr^{W^K}_{\vert K
\vert - r}(\cap_{i \in K} ker N_i) $. \\
inducing for each $(m.)$ precisely the inverse  of\\ $(-1/
(n_1^{(m_1-1)}\ldots n_n^{(m_n-1)})N_1^{(m_1-2)}\ldots
N_n^{(m_n-2)}: P((m_1, \ldots, m_n))
  \to Gr^{W^n}_{-m_n+2 } \cdots Gr^{W^1}_{-m_1+2 }(\cap_{i \in K} ker N_i)$.}\\

\n 1)
 We start the proof for $n = 1$, that is one dimensional
nilpotent orbit $(L,N)$. Then we need to prove that the
differential is the inverse of \\
\noindent \centerline {$-\frac{1}{n^r}N^{r-1}: Gr^W_{r-1}L/NL
\simeq   P(r-1)L \to Gr^W_{-(r-1)}ker N$}

\n which can be checked on the diagram
$$  \begin{array}{ccccccc}
Gr^W_{-i+2} L &&Gr^W_{-i+4} L &  \cdots &Gr^W_{i-2}
L&&Gr^W_{i} L \\
N \downarrow & \eta \searrow&N \downarrow& \cdots &N \downarrow &
\eta \searrow&N \downarrow\\
Gr^W_{-i} L&&Gr^W_{-i+2} L&\cdots&Gr^W_{i-4} L &&Gr^W_{i-2} L
\end{array}$$
where $\eta = -nId$. For $a \in Gr^W_{i-2} L$ primitive, the
element $\sum_{0 \leq  j \leq i-2} (1/n)^{i-1-j}N^{i-2-j}(a) \in
\oplus Gr^W_{i-2+2j}L$ is a cohomology class modulo the complex
 $C^{\{0\}}_{i-1} L$ inducing the cohomology class $((1/n)^{i-1}N^{i-2}(a)$ in
 $C^{\{0\}}_{1-i} L$ whose image by $ \eta$ is $-a$ the original
 primitive element up to sign.\\
 \noindent 2) In general we notice that $(N_K)^{r-2\vert K
\vert}$ decomposes on $ P(m.) L$ as the product
\\ $(n_n N_n)^{m_n-2}\circ \cdots \circ (n_1 N_1)^{m_1-2}$.\\ We
use this relation to give an
inductive proof on the number of endomorphisms $N_i$.\\
 The cohomology
of an elementary complex $K(m_1,\cdots, m_n)$ is isomorphic to the
cohomology of a complex with the unique endomorphism $N_n$ acting
on $L_n = \cap_{i \notin J(m.), i \neq n} ker N_i: L/(\Sigma_{i
\in J(m.)-\{n\}} N_i L)$ with $L_n$ depending on $m_i $ for $i
\neq n$, then the diagram is similar to the case of one variable
until we reach $-m_n +2$ for which
 the morphism $N_n^{m_n-2}$ is needed. The morphisms $N_i^{m_i-2}$
will appear  inductively with the variable $i$. So we can deduce
in general:\\
$ N_1^{m_1-2} \cdots N_n^{m_n-2}: Gr^{W^n}_{m_n } \cdots
 Gr^{W^1}_{m_1-2 } P(m.)L \simeq
Gr^{W^n}_{-m_n+2 } \cdots Gr^{W^1}_{-m_1 +2} (\cap_{i \in K} ker
N_i)$

\noindent the sum over
 $\{ ( m_1 \geq 2, \cdots, m_n \geq 2 ) :
 \Sigma_{i \in K} m_i = i +2d -1 + \vert  K \vert  $ \} induces an isomorphism:

$ \gamma: Gr^{W^K}_{i +2d -1 - \vert  K \vert}[L/ (\Sigma_{i \in
K} N_i L)] \rightarrow  Gr^{W^K}_{-(i-1)-2d +\vert  K \vert
}[(\cap_{i \in K} ( ker N_i: L \rightarrow L ] $.

 \smallskip
\noindent {\it Example}. Consider $L,N_1,N_2$ in dimension $2$, $K
= \{1,2\}$ the origin in $\C^2$, $A^K_5$, $ a \in Gr^{W^K}_2 L $,
$(m.) = (m_2= 4, m_1 = 2) $ with the following conventions for
differentials
 : the restriction from $s. = K \supset  \{1\} $ to $K$ is $-I$ ($I$ is Identity),
 the restriction from $s. = K \supset  \{2\} $ to $K$ is $I$,
  the differentials on $C^K_rL$ are ($-N_1$ on $dz_2$,$ N_i$
otherwise and $\eta$ is $n_1$  on $dz_2$, $-n_i$ otherwise.\\
An element $a$ of $A^K_5$ is written as the sum of various
components of the underlying groups in  $C^K_{-4} L \oplus
C^K_{-2}L \oplus C^K_{0} L \oplus C^K_{2}L \oplus C^K_{4}L$,  that
is  $a = \sum a(dz_J,s.,r)$
 where $dz_J$ stands for
$\wedge_{i\in J} dz_i$, $s.$ is  $\{1,2\} \supset \{1\}$, or
$\{1,2\} \supset \{2\}$ or $K$ for $\{1,2\}$ and $r$ for $a(dz_J,
s.,r) \in C^K_r L $ ($r = -4, -2, 0, 2, 4$).
 Still we need  to specify for an
element $b \in Gr^{W^K}_r L \simeq \oplus_{m_1+m_2 = r}
Gr^{W^1}_{m_1}Gr^{W^2}_{m_2} L $ its components $b =
\Sigma_{m_1+m_2 = r}b(m_1,m_2)$. A bi-primitive element $a \in
P(2,4) \subset Gr^{W^1}_0Gr^{W^2}_2L \subset Gr^{W^K}_2 L, a  = a
( dz_K,(s. = K \supset \{1\}, 4)$ in $C^K_{4}L$ (hence $N_1 a =0,
N_2^3 a = 0$ ) defines the following cohomology class $ \beta (a)
\in A^K_5/C^K_{4}L$
modulo $C^K_{4}L$: \\
$ [N_2^2 a (m_2=-2,m_1=0)(dz_{\emptyset}, s. = K, -4)\in
Gr^{W^K}_{-2} L $ in $C^K_{-4}L$ , $( n_2 N_2
a((0,0)(\emptyset, K, -2)\in Gr^{W^K}_0 L$ in $C^K_{-2}L$,\\
 $ n_1
N_2^2 (a) (0,-2)(dz_1 , s. = K \supset \{2\}, -2)\in Gr^{W^K}_{-2}
L$ in $C^K_{-2}L$, $  n_1 n_2 N_2 (a) (0,0) (dz_1, (s. = K \supset
\{1\}, 0)\in Gr^{W^K}_2 L$ in $C^K_{0}L (s. = K \supset \{1\})$,
$( n_2)^2 a((2,0)(dz_{\emptyset}, s. = K, 0)\in Gr^{W^K}_{2} L $
in $C^K_{0}L$, $(n_2)^2 n_1 a(2,0) (dz_1,s. =(K \supset \{2\}),2)$
in $C^K_{2}L(s.= K \supset \{2\}), -(n_2)^3 a(2,0) (dz_2,s.
=(K\supset \{2\}),2 )$ in $C^K_{2}L(s. = K\supset \{2\} )]$
 which ends the components of $\beta (a)$,\\
whose image by $\eta$ has two components $-(n_2)^3 n_1 a((2,0)
(dz_K,s. =(K \supset \{2\}),4)$ in $C^K_{4}L(s.=K,2) $ and $ -
(n_2)^3 n_1 a((2,0) (dz_K,s. =(K \supset \{1\}),4)$ in
$C^K_{4}L(s.= (K \supset \{1\}) $ which represents the cycle $- a$
in  the cohomology as represented by the complex
defined by $T(r)$ with $r = 4$.\\
Notice that the conditions for $s. =(K\supset \{1\}), r = 2$ are ,
$ W_1(L)\leq 2, W_1(L dz_1)\leq 0, W_1(L d z_2) \leq 2$ are
satisfied by $ a (dz_2)$ while the conditions for $s. =(K\supset
\{2\}), r = 2), r = 2$ are  $ W_2\leq 2, W_2(dz_1)\leq 2, W_2(d
z_2) \leq 0$ are not satisfied since $ W_2(a d z_2) = 2$ which
forces the lifting in
 $C^K_{2}L (s. =(K\supset \{1\})$.

\smallskip
\noindent {\bf Corollary 1}. {\em For all $i> 0$, the complex
 $ A^{KM}_i \colon= s[C^{KM}_{i+2p-1} L [p], \eta]_{-i < p
\leq 0}$ is acyclic}.\\
{\it Proof.} We can easily check, as in the previous open case,
that the cohomology of $ A^{KM}_i$ is quasi-isomorphic to the
stalk  at $y$ of the intermediate extension of the local system
 on $Y_K^* $ defined by the cohomology of $ A^{K}_i$, hence it
  is quasi-isomorphic to zero
 since $ A^{K}_i \cong 0$.

 \smallskip
\noindent The iterated monodromy morphism defines an exact
sequence\\
  $0  \rightarrow ker \, \nu^i_X \rightarrow \Psi^u_f
{\cal L}_X {\buildrel {\nu^i_X } \over \longrightarrow}
 \Psi^u_f {\cal L}_X  \rightarrow 0$ where $ ker \, \nu^i_X =
 s({\Omega}^*_X (Log Y) [p], \eta)_{-i <p \leq 0}$,\\
   $0  \rightarrow ker \, \nu^i \rightarrow \Psi^u_f
{\cal L} {\buildrel {\nu^i } \over \longrightarrow}
 \Psi^u_f {\cal L} \rightarrow 0$ where $ ker \, \nu^i =
s({\Omega}^* {\cal L} [p], \eta)_{-i <p \leq 0}$.

  \smallskip
\noindent  { \bf Corollary 2}. {\em
 For all $i \geq 1, \; Gr^{\cal W}_i ker \; \nu^i =  s(Gr^{\cal
W}_{i+2p-1}{\Omega}^* L [p], \eta)_{-i <p \leq 0} $ is acyclic.}\\
\n Proof: By the decomposition (19) we have: $ Gr^{\cal W}_i ker
\; \nu^i \cong \oplus_{K \subset M} A^{KM}_i $.

  \smallskip
\noindent  { \bf Corollary 3}. {\em For all $i\geq 1, \; \nu^i :
Gr^{\cal W}_i ({\Psi}^0_M L)
 \simeq  Gr^{\cal W}_{-i} ({\Psi}^0_M L)$.}\\
The equivalence between Corollaries 2 and  3 follows from the
exact sequence\\
  $0  \rightarrow Gr^{\cal W}_i ker \; \nu^i \rightarrow
Gr^{\cal W}_{-i} ({\Psi}^0_M L)
 {\buildrel {\nu^i } \over \longrightarrow}
  Gr^{\cal W}_{-i} ({\Psi}^0_M L)  \rightarrow 0$.

\smallskip
\noindent The  theorem follows from the corollaries.
 \be{The global weighted complex  $({\Psi}^u_f ({\cal L}), {\cal W}, F)$}
{5.3} \

\n Returning to the global situation, we need to define the Hodge
filtration on $\Psi^u_f ({\cal L}_X)$. We could define
 $$\Psi^u_f ({\cal L}_X)   \cong s(i_Y^{*}
 ({\Omega}^*_X (Log Y) \otimes {\cal L}_X)[i+1], \eta)_{i\geq 0},
  \; ({\Psi}^{u}_f {\cal L}) \colon= s({\Psi}^{u}_f {\cal L}(s.))_{s.
\in S(M)}$$
 First $F$ extends to the logarithmic complex by the formula under (15) in (4.1), then $F$
extends to $\Psi^u_f ({\cal L}_X)$  via the formula
 $$ F^r ( s(i_Y^{*}({\Omega}^*_X (Log Y) \otimes {\cal L}_X)[i+1])_{i\geq 0}) =
 s(F^{r+i+1}(i_Y^{*}(\Omega^*_X (Log Y) \otimes {\cal L}_X )[i+1], \eta)_{i\geq 0}$$
The definition of the global weight filtration reduces to the
local construction at a point  $y \in Y^*_M$, using the
quasi-isomorphism ${(\Psi}^u_f \;{\cal L})_y \; \simeq
 \;{\Psi}^0 L$.

\n We suppose again ${\cal L}$ unipotent and define as above $
{\Psi}^{u}_f {\cal L} $
 which can be viewed also as \\
 $s(i_Y^{*}{\Omega}^* {\cal L} [i+1],
 \eta)_{i \geq 0}$,
then we define the weight and Hodge filtrations
$$ {\cal W}_r ({\Psi}^{u}_f {\cal L})=
 s (i_Y^*({\cal W}_{r+2i+1} {\Omega}^* {\cal L}) [i+1], \eta)_{i
\geq 0}, \; F^r ({\Psi}^{u}_f {\cal L}) = s (i_Y^{*} F^{r+i+1}
{\Omega}^* {\cal L}[i+1], \eta)_{i \geq 0}$$ The logarithm of the
monodromy ${\cal N}$ is defined on this complex as in the local
case. The filtration $W({\cal N})$ is defined on ${\Psi}^{u}_f
{\cal L}$ in the abelian category of perverse sheaves.
  \be {Theorem.}{} {\em Suppose $\cal L$ underlies
a unipotent variation of polarised Hodge structures of weight $m$,
then the graded part of the weight filtration  of the complex}

\centerline {$(\Psi^u_f {\cal L},{\cal W}[m], F)$}

\n {\em decomposes into a direct sum of intermediate extension of
polarised $ VHS$. Moreover it is a $MHC$ for X proper.\\
We  have
 $W({\cal N})= {\cal W}$; $ {\cal W}_{r+1} {\Omega}^* {\cal L}$
 is induced by $ {\cal W}_r \Psi^u_f {\cal L}$ on $Coker {\cal N}$
 and $Gr^{\cal W}_{r+1}i_Y^{*}{\Omega}^* {\cal L}$
 is the primitive part of $ Gr^{\cal W}_{r}\Psi^u_f {\cal L}$. }

\smallskip
\n The proof of this theorem reduces by definition to show that
$(Gr^{{\cal W}[m]}_r (\Psi^u_f ({\cal L}_X), F) $ decomposes which
result can be reduced to the local case where it follows from

\smallskip
\noindent { \bf Lemma}. {\em Let $i_0 $ be a positive integer
large enough to have $ Gr^{\cal W}_{j}{\Omega}^* {\cal L} = 0$ for
all $ \vert j \vert > i_0$ and  $ I(r) = \{ p \geq 0 : \vert r
\vert + 1 \leq r + 2p + 1 \leq i_0 \}$, then}
$$Gr^{\cal W}_{r}({\Psi}^0 L)_{M} =
 s(Gr^{\cal W}_{ r+2p-1}{\Omega}^* L [p], \eta)_{ p\in I(r)}; \,
   Gr^{\cal W}_{r}({\Psi}^0 L)_{ M} \stackrel {\nu^r }{\stackrel {\cong}\rightarrow}
 Gr^{\cal W}_{-r}({\Psi}^0 L)_{ M}, r \geq 0 $$
\noindent Proof:  Suppose $r > 0$, then $Gr^{\cal W}_{r}({\Psi}^0
L)_{ M} \cong s(Gr^{\cal W}_{r+2p+1}{\Omega}^* L [p], \eta)_{p \in
I(r)}$ since $ Gr^{\cal W}_{r+2p+1}{\Omega}^* L \cong 0$ for $r
\notin I(r)$, while for $r < 0$, the sum\\
 $ s(Gr^{\cal W}_{r+2p+1}{\Omega}^* L [p], \eta)_{ -1 < p < -r }
 \cong s(Gr^{\cal W}_{i}{\Omega}^* L [p], \eta)_{i \in
 [r+1, -(r+1)]}$ is acyclic.\\
 The quasi-isomorphism defined by $\nu^r$ follows since the
 formula is symmetric in $r$ and $-r$.

\smallskip
\noindent Example. $ Gr^{\cal W}_{0}({\Psi}^0 L)_{M} \cong
s(Gr^{\cal W}_{1}{\Omega}^* L [1] \stackrel {\eta }{\rightarrow}
Gr^{\cal
W}_{3}{\Omega}^* L [2]\ldots) $\\
$ Gr^{\cal W}_{1}({\Psi}^0 L)_{ M} \cong s(Gr^{\cal
W}_{2}{\Omega}^* L [1] \stackrel {\eta }{\rightarrow}
Gr^{\cal W}_{4}{\Omega}^* L [2]\ldots ) $\\
$ Gr^{\cal W}_{2}({\Psi}^0 L)_{ M} \cong s(Gr^{\cal
W}_{3}{\Omega}^* L [1] \stackrel {\eta }{\rightarrow} Gr^{\cal
W}_{5}{\Omega}^* L [2]\ldots ) $\\
For $r = i_0-1$ and $ r= i_0-2$, $ Gr^{\cal W}_{r}({\Psi}^0 L)_{
M} \cong Gr^{\cal W}_{i_0}{\Omega}^* L [1]$

\smallskip
\noindent { \bf Corollary}. {\em  The  graded part of $({\Psi}^0
L)_{ M}$ is non zero for only a finite number of indices  which
decomposes into direct sum of of intermediate extension of
polarised $ VHS$}.

\smallskip
\noindent Proof. We can introduce the notion of primitive parts
starting with $ Gr^{\cal W}_{r}({\Psi}^0 L)_{ M}$ for $r = i_0-1$
and $ r= i_0-2$ which are direct sum of of intermediate extension
of polarised $ VHS$, then using the decomposition into direct sum
of primitive parts we prove by induction that such decomposition
is
valid for all indices $r$ (there is no extension by $\eta$).\\
Then the theorem follows easily.

\smallskip
\noindent{\bf Remark}({\it dual statement}). \n  Let us return to
  $ ({\Psi}^{u}_f {\cal L} \colon= s({\Psi}^{u}_f {\cal
L}(s.))_{s. \in S(I)}$ defined with indices $p \leq 0$,
 which can be viewed locally at a point $y \in Y$  as $s({\Omega}^* {\cal L} [p],
 \eta)_{p \leq 0}$,
 the weight and Hodge filtrations  defined by the formula (18).

 \smallskip
\noindent {\bf Theorem.}({\it dual statement}) {\em Suppose $\cal
L$ underlies a unipotent variation of polarised Hodge structures
of weight $m$, then the graded part of the weight filtration  of
the complex }

\centerline {$(\Psi^u_f {\cal L},{\cal W}[m], F)$}

\n {\em decomposes into a direct sum of intermediate extension of
polarised $ VHS$. Moreover it is a $MHC$ for X proper.\\
We  have
 $W({\cal N})= {\cal W}$; $ {\cal W}_{r-1} i_Y^{*}{\Omega}^* {\cal L}$
 is induced by $ {\cal W}_{r} \Psi^u_f {\cal L}$ on $ Ker {\cal N}$
 for $r \leq 0$  and $ Gr^{\cal W}_{r-1}i_Y^{*}{\Omega}^* {\cal L}$
 is isomorphic to the primitive part of
  $ Gr^{\cal W}_{-r}\Psi^u_f {\cal L}$. }

\smallskip
\n The proof of this theorem reduces by definition to show that
$(Gr^{{\cal W}[m]}_r \Psi^u_f {\cal L}, F) $ decomposes; such
result can be reduced to the local case where it follows from

\smallskip
\noindent { \bf Lemma}. {\em Let $i_0 $ be a positive integer
large enough to have $ Gr^{\cal W}_{j}{\Omega}^* {\cal L} = 0$ for
all $ \vert j \vert > i_0$ and  $ I(r) = \{ p \leq 0 ,  - i_0 + 1
\leq r + 2p -1 \leq - \vert r \vert -1 \}$, then }
$$Gr^{\cal W}_{r}({\Psi}^0 L)_{ M} =
 s(Gr^{\cal W}_{r+2p-1}{\Omega}^* L [p], \eta)_{ p\in I(r)}, \;
   Gr^{\cal W}_{r}({\Psi}^0 L)_{ M} \stackrel {\nu^r }
   {\stackrel {\cong}\rightarrow}
 Gr^{\cal W}_{-r}({\Psi}^0 L)_{ M}, r \geq 0 $$
\noindent Proof: Suppose $r > 0$, then $Gr^{\cal W}_{r}({\Psi}^0
L)_{ M}$  is the sum of
 $ s(Gr^{\cal W}_{r+2p-1}{\Omega}^* L [p], \eta)_{ -r < p
\leq 0}\cong 0$ and $ s(Gr^{\cal W}_{r+2p-1}{\Omega}^* L [p],
\eta)_{ p \leq -r} $. In particular, $Gr^{\cal W}_{r}({\Psi}^0
L)_{ M} \cong 0$ for all $r$ such that  $ \vert r \vert \geq i_0
$. Using (21 ) and (22) we get a decomposition into $ \Psi^{KM}_r
L = s( C^{KM}_{r+2p-1} L [p], \eta )_{p \in I(r)}$ where
$C^{KM}_{r+2p-1} L$ is the fiber of an intermediate extension of a
local system defined by $C^{K}_{r+2p-1} L$ for $r+2p-1 < 0 $ on
$Y^*_K$. \\
Then the proof is similar to the previous case in the dual
definition of $\Psi^u_f {\cal L}$.

\medskip
\n ACKNOWLEDGMENTS. I would like to acknowledge the valuable
discussions I  had  with J.L. Verdier at the time I started to
work on this subject and to thank J.P. Serre for accepting the
note [14].
 
\hfill Fouad El Zein

 \quad  \hfill  {Projet de G\'eom\'etrie et
Dynamique du CNRS, Institut de Math\'ematiques 
 Jussieu.}

\end{document}